# Extensions to Generalized Disjunctive Programming: Hierarchical Structures and First-order Logic

Hector D. Perez, Ignacio E. Grossmann*
*Department of Chemical Engineering, Carnegie Mellon University, Pittsburgh, PA, 15213, USA*

Abstract

Optimization problems with discrete-continuous decisions are traditionally modeled in algebraic form via (non)linear mixed-integer programming. A more systematic approach to modeling such systems is to use Generalized Disjunctive Programming (GDP), which extends the Disjunctive Programming paradigm proposed by Egon Balas to allow modeling systems from a logic-based level of abstraction that captures the fundamental rules governing such systems via algebraic constraints and logic. Although GDP provides a more general way of modeling systems, it warrants further generalization to encompass systems presenting a hierarchical structure. This work extends the GDP literature to address three major alternatives for modeling and solving systems with nested (hierarchical) disjunctions: *explicit nested disjunctions*, *equivalent single-level disjunctions*, and *flattening via basic steps*. We also provide theoretical proofs on the relaxation tightness of such alternatives, showing that explicitly modeling nested disjunctions is superior to the traditional approach discussed in literature for dealing with nested disjunctions.

1. Introduction

Discrete-continuous optimization is one of the main modeling approaches to address design, planning, and scheduling problems in Process Systems Engineering (PSE) (Grossmann, 2012). Raman and Grossmann (1994) present a powerful modeling paradigm that extends the work by Balas (1985) on disjunctive programming. This new paradigm, called Generalized Disjunctive Programming (GDP), has been further developed by others in the PSE community over the years to account for additional features, such as nonlinearities and nonconvexities in the problems encountered (Grossmann & Trespalacios, 2013). GDP relies on the intersection of disjunctions of convex sets to model the feasible space. Boolean variables are used as indicator variables for each convex set, meaning that if *True*, the constraints in the disjunct are enforced. Logic constraints are also included to describe the relationships between the Boolean indicator variables via propositional logic and constraint programming logic.

GDP is a powerful modeling abstraction for optimization problems for two main reasons. Firstly, modeling systems from the basis of their underlying logical relationships speeds up the development of optimization models by making them easier to interpret, and reducing the likelihood of modeling errors due to logical fallacies. Secondly, GDP makes available a broad array of solution methods, ranging from mixed-integer reformulations to logic-based search methods (Chen et al., 2022).

The present work extends the GDP theory to allow modeling hierarchical systems, which are commonly encountered in PSE, and more particularly in Enterprise-Wide Optimization (EWO) (Grossmann, 2012; van

den Heever & Grossmann, 1999), and flowsheet superstructure optimization (Türkay & Grossmann, 1996a). Hierarchical systems involve multiple level of decision making, which can be concisely modelled via nested disjunctions. However, traditional GDP does not consider such formulations. Existing GDP literature suggests reformulating nested disjunctions into equivalent single-level disjunctions (Vecchietti & Grossmann, 2000). Such an approach requires introducing additional Boolean variables and logical propositions. An alternate approach is used in the work by van den Heever and Grossmann (1999), in which a direct or inside-out reformulation to MI(N)LP is performed. We formalize these two approaches and provide theoretical proofs on the tightness of their continuous relaxations. We also present a third approach that produces a single-level GDP by applying basic steps to obtain the disjunctive normal form (DNF) of the nested disjunction. The model tightness and computational performance of the different approaches are compared. A series of examples are used to show the modeling and computational advantages obtained by explicitly modeling nested disjunctions.

The paper is organized as follows, Section 2 provides a background on the GDP modeling paradigm. Section 3 extends this formulation to account for hierarchical systems, and discusses the alternatives for modeling such systems. The equivalent mixed-integer programming reformulations for these alternatives are presented, along with two theorems on the tightness of the resulting models. Section 4 provides several numerical use cases for hierarchical GDPs. Section 5 presents concluding remarks.

## 2. Background: Generalized Disjunctive Programming (GDP)

The classical GDP formulation is given below (*GDP*), where $x$ is the set of continuous variables (bounded between $x^{LB}$ and $x^{UB}$), $f(x)$ is the objective function, $r(x) \leq 0$ is the set of global constraints, $g_{ij}(x) \leq 0$ is the set of constraints applied when the indicator Boolean $Y_{ij}$ is *True* for disjunct $j$ in disjunction $i$. $f(x)$, $r(x)$, and $g_{ij}(x)$ are assumed to be continuous and differentiable over $x$. $\Omega(Y)$ defines the set of logic constraints, which are described via propositional logic on a subset of Boolean variables. These constraints describe the relations between the Boolean variables via clauses that contain with one or more of the following logic operators: AND (∧), OR (∨), implication (⇒), equivalence (⇔), and negation (¬). The set of logic constraints may also include cardinality clauses of the form *choose exactly* (or *at least* or *at most*) $m$ Boolean variables from a subset of Booleans to be *True* (Yan & Hooker, 1999). We leverage predicate logic to extend the notation used by Yan and Hooker for cardinality clauses by defining the following predicates: $\Xi(m, Y_s \ \forall s \in S)$ enforces *exactly* $m$ of the Boolean variables $Y_s$ are *True*, $\Lambda(m, Y_s \ \forall s \in S)$ enforces *at least* $m$ of the variables are *True*, and $\Gamma(m, Y_s \ \forall s \in S)$ enforces *at most* $m$ *are True*. GDP models typically include a cardinality clause to enforce that exactly 1 disjunct is chosen in each disjunction, i.e., $\Xi(1, Y_{ij} \ \forall j \in J_i) \ \forall i \in I$. The GDP literature often uses the exclusive OR operator, $\underline{\vee}$, to define this constraint. However, such an operator is only correct for proper disjunctions (those with non-overlapping disjuncts). Thus, to avoid any ambiguity, we use the predicate logic notation, $\Xi(1, Y)$, here instead.

$$\min z = f(x) \qquad \text{(GDP)}$$
$$s.t. \quad r(x) \leq 0$$
$$\bigvee_{j \in J_i} \begin{bmatrix} Y_{ij} \\ g_{ij}(x) \leq 0 \end{bmatrix} \qquad \forall i \in I$$
$$\Xi(1, Y_{ij} \ \forall j \in J_i) \qquad \forall i \in I$$

$$\Omega(Y)$$
$$x^{LB} \leq x \leq x^{UB}$$
$$x \in \mathbb{R}^n$$
$$Y_{ij} \in \{True, False\} \qquad \forall i \in I, j \in J_i$$

To illustrate the elements of a GDP model, consider the model below (*GDP-example*). The projection of this model on the $x_1, x_2$-plane is given in **Figure 2.1**, where the quadratic objective function is shown in the colored contours, the global constraints are given by the region under the black curves (one linear and the other nonlinear), and the three disjuncts given in the colored rectangles. The feasible space of such a system is given by the disjoint regions in the orange, blue, and green rectangles that satisfy the global constraints.

$$\min z = \frac{1}{2}(x_1 - 2)^2 + \frac{3}{2}(x_2 - 3)^2 \qquad \text{(GDP-example)}$$

$$s.t. \quad \frac{1}{10}x_1^2 + x_2 \leq 3$$

$$2x_1 + x_2 \leq 10$$

$$\begin{bmatrix} Y_1 \\ 0 \leq x_1 \leq 1 \\ 2 \leq x_2 \leq 3 \end{bmatrix} \vee \begin{bmatrix} Y_2 \\ \frac{3}{2} \leq x_1 \leq 1 \\ 0 \leq x_2 \leq 1 \end{bmatrix} \vee \begin{bmatrix} Y_3 \\ \frac{9}{4} \leq x_1 \leq \frac{15}{4} \\ 1 \leq x_2 \leq \frac{11}{5} \end{bmatrix}$$

$$\Xi(1, Y_i \ \forall i \in \{1,2,3\})$$
$$0 \leq x_1, x_2 \leq 5$$
$$x_1, x_2 \in \mathbb{R}^1$$
$$Y_i \in \{True, False\} \qquad \forall i \in \{1,2,3\}$$

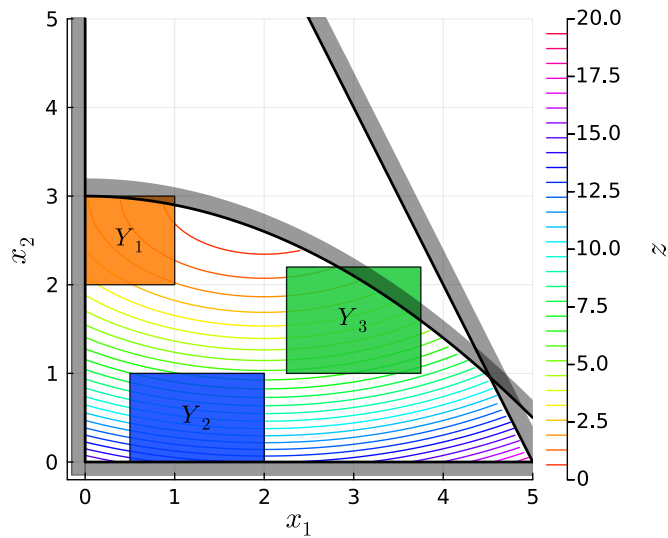

**Figure 2.1.** Sample GDP graphical representation for *GDP-example* model.

One of the main advantages of modeling discrete-continuous problems using GDP is the collection of methods that are available for optimizing such systems. These include, 1) reformulating to mixed-integer (non)linear models (MI(N)LP) via either Big-M (Trespalacios & Grossmann, 2015) or Hull reformulations (Agarwal, 2015; Bernal & Grossmann, 2021; Furman et al., 2020; Grossmann & Lee, 2003), 2) logic-based decomposition methods such as Logic-based Outer Approximation (LOA) (Türkay & Grossmann, 1996b), 3) disjunctive branch-and-bound (Lee & Grossmann, 2000), 4) basic steps (Ruiz & Grossmann, 2012), and 5) hybrid cutting planes (Sawaya & Grossmann, 2005; Trespalacios & Grossmann, 2016). The reader is referred to the above references for a detailed understanding of each of these solution methods.

3. Extended formulation for multi-level hierarchies

Decision hierarchies are present in most decision-making applications. These include for instance supply chain and enterprise-wide optimization, where different levels of decision-making exist depending on the time scales considered: planning (months/years), scheduling (hours/days), and control (seconds/minutes). According to Brunaud and Grossmann (2017), integrating different decision levels enables better coordination and communication between functional areas, which increases agility in response to disturbances and makes it possible to attain benefits for the company that are not possible with a siloed approach. **Figure 3.1** illustrates the notion of the synergistic benefits that can be obtained by an integrated approach, rather than siloed or aggregated approaches. Accounting for the relationships between different levels of decision-making can aid in finding the true optimum, which differs from that of the aggregated model (i.e., the model obtained by summing the siloed costs). Integrated approaches to hierarchical decision-making systems have been addressed in the literature. Some examples of these integrations are the integration between design and planning (operational and expansion) (van den Heever & Grossmann, 1999), planning and scheduling (Maravelias & Sung, 2009), and scheduling and control (Muñoz et al., 2011; Sokoler et al., 2017). The following subsections formalize how GDP can be used to model hierarchical systems, along with theoretical proofs on the differences between the approaches.

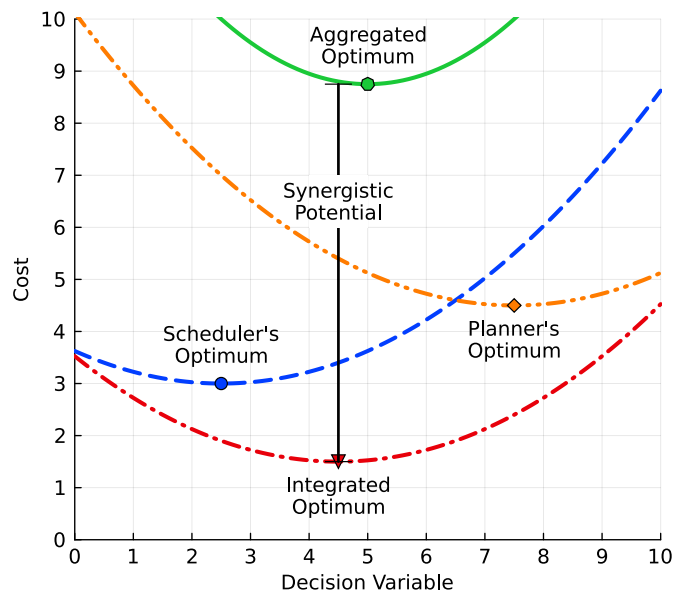

**Figure 3.1.** Illustartion of the different optimas for siloed, aggregated, and integrated approaches.

### 3.1. Hierarchical GDP

We propose extending the GDP paradigm to include multi-level decisions by means of nested disjunctions. Although the notion of nesting disjunctions to represent hierarchical decisions is not new, the limitations in the traditional GDP notation have made it difficult to exploit the benefits of using such structures. One of the first references to nested disjunctions is found in the work by Vecchietti and Grossmann (2000), which describes the transformations required to conform to the current GDP notation. It is interesting to note that several works have relied on the nested GDP representation due to its compact representation. In one of these (Rodriguez & Vecchietti, 2009), the following statement is made,

> *"Although the expressiveness of the hierarchical decisions by means of nested disjunctions, they cannot be implemented directly. These disjunctions must be transformed into GDP form. For that purpose, the disjunctions…must be rewritten as single disjunctions, and some additional constraints must also be included in the model."*

Therefore, from a model development point of view, the use of disjunction nesting is shown to add value. However, its implementation has often required breaking the explicit hierarchical structure. An exception is the work by van den Heever and Grossmann (1999), which does not transform the nested GDP into a logically equivalent single-level GDP, but rather suggests performing the Hull reformulation on the inner disjunction and then reformulating the outer disjunction. We now build upon this concept to formally extend the GDP notation for hierarchical systems that generalizes to multi-disjunct disjunctions, rather than the *on/off* disjunctions used by van den Heever and Grossmann (1999). We also provide both theoretical proofs on the advantages of modeling system hierarchies via nested disjunctions, and highlight the computational performance gains obtained using this explicit notation.

The proposed extension to the classical GDP notation for hierarchical systems is given below for a **2-Level Nested GDP** (*2L-GDP*), where the upper-level decisions, $Y$, enforce the constraints $g(x) \leq 0$ and the nested decisions, $W$, which have constraints $h(x) \leq 0$. Here the cardinality clause of selecting exactly one disjunct from the upper-level decisions, $Y$, is expressed explicitly, along with a new set of cardinality rules that enforce selecting exactly one of the lower-level decisions, $W$, if and only if the upper-level decision has been selected, and selecting no lower-level decisions when the upper-level decision is not selected. This constraint is expressed as the conjunction of two cardinality rules: $\left[Y_{ij} \Rightarrow \Xi\left(1, W_{ijkl} \ \forall l \in L_{ijk}\right)\right] \wedge \left[\neg Y_{ij} \Rightarrow \Xi\left(0, W_{ijkl} \ \forall l \in L_{ijk}\right)\right] \forall i \in I, j \in J_i, k \in K_{ij}$. In the GDP literature, this constraint has been traditionally written as $Y_{ij} \Leftrightarrow \underline{\vee}_{l \in L_{ijk}} W_{ijkl} \ \forall i \in I, j \in J_i, k \in K_{ij}$. However, such a logic proposition is incomplete because it would allow the following to occur: $Y_{ij} = False$ and $W_{ijkl} = True$ for more than 1 index $l \in L_{ijk}$ (i.e., *False* $\Leftrightarrow$ (*True* $\underline{\vee}$ *True*) is valid because the exclusive OR makes the right-hand side *False*). If all disjunctions are proper, then this will not occur. However, since there can be a disjunction with overlapping disjuncts, the cardinality rule $\Gamma\left(1, W_{ijkl} \ \forall l \in L_{ijk}\right) \forall i \in I, j \in J_i, k \in K_{ij}$ would need to be added to such a system to ensure that no more than 1 literal, $W_{ijkl}$, is set to *True*. A more compact form would be to use the predicate constraint, $\Xi\left(t(Y_{ij}), W_{ijkl} \ \forall l \in L_{ijk}\right)$, where $t(\cdot)$ is a unary function that maps a Boolean variable to its binary counterpart (i.e., $t(True) = 1$ and $t(False) = 0$). For simplicity, we make a slight abuse of notation, dropping the mapping function $t(\cdot)$, and using the expression $\Xi\left(Y_{ij}, W_{ijkl} \ \forall l \in L_{ijk}\right)$ instead.

$$\min z = f(x) \tag{2L-GDP}$$

$$s.t. \quad r(x) \leq 0$$

$$\bigvee_{j \in J_i} \left[ \begin{array}{c} Y_{ij} \\ g_{ij}(x) \leq 0 \\ \bigvee_{l \in L_{ijk}} \left[ \begin{array}{c} W_{ijkl} \\ h_{ijkl}(x) \leq 0 \end{array} \right] \forall k \in K_{ij} \end{array} \right] \quad \forall i \in I$$

$$\Xi(1, Y_{ij} \; \forall j \in J_i) \quad \forall i \in I$$
$$\Xi(Y_{ij}, W_{ijkl} \; \forall l \in L_{ijk}) \quad \forall i \in I, j \in J_i, k \in K_{ij}$$
$$\Omega(Y, W)$$
$$x^{LB} \leq x \leq x^{UB}$$
$$x \in \mathbb{R}^n$$
$$Y_{ij} \in \{True, False\} \quad \forall i \in I, j \in J_i$$
$$W_{ijkl} \in \{True, False\} \quad \forall i \in I, j \in J_i, k \in K_{ij}, l \in L_{ijk}$$

This model can be generalized to a **Multi-Level Nested GDP** (*ML-GDP*) with $n$ levels, where the superscript on the Boolean variables, constraints, and sets indicates the level $k \in \{1, ..., n\}$ of the hierarchy that these belong to.

$$\min z = f(x) \tag{ML-GDP}$$
$$s.t. \quad r(x) \leq 0$$

$$\bigvee_{j_1 \in J_{i_1}^{(1)}} \left[ \begin{array}{c} Y_{i_1 j_1}^{(1)} \\ g_{i_1 j_1}^{(1)}(x) \leq 0 \\ \bigvee_{j_2 \in J_{i_1 j_1 i_2}^{(2)}} \left[ \begin{array}{c} Y_{i_1 j_1 i_2 j_2}^{(2)} \\ g_{i_1 j_1 i_2 j_2}^{(2)}(x) \leq 0 \\ \vdots \\ \bigvee_{j_n \in J_{i_1 j_1 \dots i_n}^{(n)}} \left[ \begin{array}{c} Y_{i_1 j_1 \dots i_n j_n}^{(n)} \\ g_{i_1 j_1 \dots i_n j_n}^{(n)}(x) \leq 0 \end{array} \right] \forall i_n \in I_{i_1 j_1 \dots i_{n-1} j_{n-1}}^{(n)} \end{array} \right] \forall i_2 \in I_{i_1 j_1}^{(2)} \end{array} \right]$$
$$\forall i_1 \in I^{(1)}$$

$$\Xi(1, Y_{i_1 j_1}^1 \; \forall j_1 \in J_{i_1}) \quad \forall i_1 \in I^{(1)}$$

$$\Xi\left(Y_{i_1 j_1 \dots i_{k-1} j_{k-1}}^{(k-1)}, Y_{i_1 j_1 \dots i_k j_k}^{(k)} \; \forall j_k \in J_{i_1 j_1 \dots i_k}^{(k)}\right)$$
$$\forall k \in \{2, ..., n\}, i_1 \in I^{(1)}, j_1 \in J_{i_1}^{(1)}, ..., i_k \in I_{i_1 j_1 \dots i_{k-1} j_{k-1}}^{(k)}, j_k \in J_{i_1 j_1 \dots i_k}^{(k)}$$

$$\Omega(Y^{(1)}, ..., Y^{(n)})$$
$$x^{LB} \leq x \leq x^{UB}$$
$$x \in \mathbb{R}^n$$
$$Y_{i_1 j_1 \dots i_k j_k}^{(k)} \in \{True, False\}$$
$$\forall k \in \{1, ..., n\}, i_1 \in I^{(1)}, j_1 \in J_{i_1}^{(1)}, ..., i_k \in I_{i_1 j_1 \dots i_{k-1} j_{k-1}}^{(k)}, j_k \in J_{i_1 j_1 \dots i_k}^{(k)}$$

## 3.2. Equivalent Single-level GDP

Previous references to GDP with nested disjunctions in literature have proposed transforming the 2L-GDP model into the **Equivalent Single-Level GDP** (*2E-GDP*) given below (Grossmann & Trespalacios, 2013; Vecchietti & Grossmann, 2000). Here, the nested disjunction is extracted and a dummy or "slack" disjunct is added to preserve feasibility. Thus, if none of the nested disjuncts is selected, the slack disjunct is selected, which contains the entire feasible set for $x$. The exclusive cardinality rule on the inner Boolean variables, $W$, is also augmented to include the slack Boolean variable, $W_{ijk0}$. This slack variable is, however, not included in the linking logic constraint for the upper and lower-level decisions. This ensures that the nested decisions are only selected if their master Boolean is *True*. This method for transforming a nested disjunction can also be applied to the multi-level system *ML-GDP*.

$$\min z = f(x) \tag{2E-GDP}$$
$$s.t. \quad r(x) \leq 0$$

$$\bigvee_{j \in J_i} \begin{bmatrix} Y_{ij} \\ g_{ij}(x) \leq 0 \end{bmatrix} \qquad \forall i \in I$$

$$\left( \bigvee_{l \in L_{ijk}} \begin{bmatrix} W_{ijkl} \\ h_{ijkl}(x) \leq 0 \end{bmatrix} \right) \bigvee \begin{bmatrix} W_{ijk0} \\ x^{LB} \leq x \leq x^{UB} \end{bmatrix} \qquad \forall i \in I, j \in J_i, k \in K_{ij}$$

$$\Xi(1, Y_{ij} \; \forall j \in J_i) \qquad \forall i \in I$$
$$\Xi(1, W_{ijkl} \; \forall l \in L_{ijk} \cup \{0\}) \qquad \forall i \in I, j \in J_i, k \in K_{ij}$$
$$\Xi(Y_{ij}, W_{ijkl} \; \forall l \in L_{ijk}) \qquad \forall i \in I, j \in J_i, k \in K_{ij}$$
$$\Omega(Y, W)$$
$$x^{LB} \leq x \leq x^{UB}$$
$$x \in \mathbb{R}^n$$
$$Y_{ij} \in \{True, False\} \qquad \forall i \in I, j \in J_i$$
$$W_{ijkl} \in \{True, False\} \qquad \forall i \in I, j \in J_i, k \in K_{ij}, l \in L_{ijk}$$

Although the above formulation, allows modeling hierarchical systems in the standard GDP notation, it has two major drawbacks: 1) the explicit hierarchical structure is lost, and 2) although the Equivalent Single-Level GDP model is logically equivalent to the Nested GDP model, it requires introducing additional disjuncts and Boolean variables. Introducing "slack" disjuncts and "slack" Boolean variables results in models whose continuous relaxations are less tight, as described in the next section.

## 3.3. Tightness of Continuous Relaxations

The following two theorems and their associated proofs establish the advantages of modeling multi-level decisions problems via Nested GDP, rather than the Equivalent Single-Level GDP approach. The advantages are shown by discussing the tightness of the continuous relaxations of both the Hull reformulation (HR) and Big-M reformulation (BM) of these two GDP models.

*Theorem 1.* Let *rML-GDP-HR* denote the continuous relaxation of the mixed-integer program (MIP) obtained from a Multi-Level Nested GDP via the Hull reformulation, and let *rME-GDP-HR* denote the continuous relaxation of the MIP obtained from its respective Equivalent Single-Level GDP representation via the Hull reformulation. The feasible space of the former is contained within the feasible space of the latter, namely, *rML-GDP-HR* ⊆ *rME-GDP-HR*.

*Proof.* Without loss of generality, the above theorem is proved by establishing that the Hull reformulation of the 2-Level Nested GDP model (*r2L-GDP-HR*) is contained in the Hull reformulation of its Equivalent Single-Level GDP representation (*r2E-GDP-HR*):

$$r2L\text{-}GDP\text{-}HR \subseteq r2E\text{-}GDP\text{-}HR$$

The Hull reformulation for *2L-GDP* is given below, where the continuous variable $x$ is disaggregated in each disjunct ($x$ is disaggregated into $u_{ij}$ for each upper-level disjunct, and $u_{ij}$ is disaggregated into $v_{ijkl}$ for each lower-level disjunct) and the Boolean variables are replaced by their corresponding binary variable ($Y$ becomes $y$, and $W$ becomes $w$). $A$ and $B$ are matrices of scalars, and $c$ is a vector of scalars. These are used to map the logic constraints into their algebraic counterparts obtained after converting the logic propositions into conjunctive normal form (CNF) and transforming each clause into its equivalent algebraic constraint (Williams, 1985).

$$\min z = f(x) \qquad \text{(2L-GDP-HR)}$$

$$s.t. \quad r(x) \leq 0$$

**Variable Disaggregation:**

$$\sum_{j \in J_i} u_{ij} = x \qquad \forall i \in I$$

$$\sum_{l \in L_{ijk}} v_{ijkl} = u_{ij} \qquad \forall i \in I, j \in J_i, k \in K_{ij}, L_{ijk} \neq \emptyset$$

$$x^{LB} \cdot y_{ij} \leq u_{ij} \leq x^{UB} \cdot y_{ij} \qquad \forall i \in I, j \in J_i$$

$$x^{LB} \cdot w_{ijkl} \leq v_{ijkl} \leq x^{UB} \cdot w_{ijkl} \qquad \forall i \in I, j \in J_i, k \in K_{ij}, l \in L_{ijk}$$

**Perspective Functions:**

$$y_{ij} \cdot g_{ij}\left(\frac{u_{ij}}{y_{ij}}\right) \leq 0 \qquad \forall i \in I, j \in J_i$$

$$w_{ijkl} \cdot h_{ijkl}\left(\frac{v_{ijkl}}{w_{ijkl}}\right) \leq 0 \qquad \forall i \in I, j \in J_i, k \in K_{ij}, l \in L_{ijk}$$

**Cardinality Rules:**

$$\sum_{j \in J_i} y_{ij} = 1 \qquad \forall i \in I$$

$$y_{ij} = \sum_{l \in L_{ijk}} w_{ijkl} \qquad \forall i \in I, j \in J_i, k \in K_{ij}, L_{ijk} \neq \emptyset$$

| Other Logic Constraints | $Ay + Bw \leq c$ | |
|---|---|---|

| Variable Domains | $x^{LB} \leq x \leq x^{UB}$ | |
|---|---|---|
| | $\min(0, x^{LB}) \leq u_{ij} \leq \max(0, x^{UB})$ | $\forall i \in I, j \in J_i$ |
| | $\min(0, x^{LB}) \leq v_{ijkl} \leq \max(0, x^{UB})$ | $\forall i \in I, j \in J_i, k \in K_{ij}, l \in L_{ijk}$ |
| | $x \in \mathbb{R}^n$ | |
| | $u_{ij} \in \mathbb{R}^n$ | $\forall i \in I, j \in J_i$ |
| | $v_{ijkl} \in \mathbb{R}^n$ | $\forall i \in I, j \in J_i, k \in K_{ij}, l \in L_{ijk}$ |
| | $y_{ij} \in \{0,1\}$ | $\forall i \in I, j \in J_i$ |
| | $w_{ijkl} \in \{0,1\}$ | $\forall i \in I, j \in J_i, k \in K_{ij}, l \in L_{ijk}$ |

The Hull reformulation for *2E-GDP* is given below, where $x$ is disaggregated into $u_{ij}$ for the upper-level disjunctions, and is also disaggregated into $v_{ijkl}$ for the lower-level disjuncts, which are extracted when transforming the model into an Equivalent Single-Level GDP.

$$\min z = f(x) \quad \text{(2L-GDP-HR)}$$
$$s.t. \quad r(x) \leq 0$$

| Variable Disaggregation | $\sum_{j \in J_i} u_{ij} = x$ | $\forall i \in I$ |
|---|---|---|
| | $v_{ijk0} + \sum_{l \in L_{ijk}} v_{ijkl} = x$ | $\forall i \in I, j \in J_i, k \in K_{ij}, L_{ijk} \neq \emptyset$ |
| | $x^{LB} \cdot y_{ij} \leq u_{ij} \leq x^{UB} \cdot y_{ij}$ | $\forall i \in I, j \in J_i$ |
| | $x^{LB} \cdot w_{ijkl} \leq v_{ijkl} \leq x^{UB} \cdot w_{ijkl}$ | $\forall i \in I, j \in J_i, k \in K_{ij}, l \in L_{ijk} \cup \{0\}, L_{ijk} \neq \emptyset$ |

| Perspective Functions | $y_{ij} \cdot g_{ij}\left(\dfrac{u_{ij}}{y_{ij}}\right) \leq 0$ | $\forall i \in I, j \in J_i$ |
|---|---|---|
| | $w_{ijkl} \cdot h_{ijkl}\left(\dfrac{v_{ijkl}}{w_{ijkl}}\right) \leq 0$ | $\forall i \in I, j \in J_i, k \in K_{ij}, l \in L_{ijk}$ |

| Cardinality Rules | $\sum_{j \in J_i} y_{ij} = 1$ | $\forall i \in I$ |
|---|---|---|
| | $y_{ij} = \sum_{l \in L_{ijk}} w_{ijkl}$ | $\forall i \in I, j \in J_i, k \in K_{ij}, L_{ijk} \neq \emptyset$ |
| | $w_{ijk0} + \sum_{l \in L_{ijk}} w_{ijkl} = 1$ | $\forall i \in I, j \in J_i, k \in K_{ij}, L_{ijk} \neq \emptyset$ |

| | |
|---|---|
| **Other Logic Constraints** | $Ay + Bw \leq c$ |
| **Variable Domains** | $x^{LB} \leq x \leq x^{UB}$ <br> $\min(0, x^{LB}) \leq u_{ij} \leq \max(0, x^{UB})$ $\qquad \forall i \in I, j \in J_i$ <br> $\min(0, x^{LB}) \leq v_{ijkl} \leq \max(0, x^{UB})$ $\quad \forall i \in I, j \in J_i, k \in K_{ij}, l \in L_{ijk} \cup \{0\}, L_{ijk} \neq \emptyset$ <br> $x \in \mathbb{R}^n$ <br> $u_{ij} \in \mathbb{R}^n$ $\qquad \forall i \in I, j \in J_i$ <br> $v_{ijkl} \in \mathbb{R}^n$ $\qquad \forall i \in I, j \in J_i, k \in K_{ij}, l \in L_{ijk} \cup \{0\}, L_{ijk} \neq \emptyset$ <br> $y_{ij} \in \{0,1\}$ $\qquad \forall i \in I, j \in J_i$ <br> $w_{ijkl} \in \{0,1\}$ $\qquad \forall i \in I, j \in J_i, k \in K_{ij}, l \in L_{ijk} \cup \{0\}, L_{ijk} \neq \emptyset$ |

The difference between *2L-GDP-HR* and *2E-GDP-HR* is in the highlighted constraints in the variable disaggregation and cardinality rules sections. The proof for the Hull reformulation case is given by applying Fourier-Motzkin elimination (Dantzig, 1972) to eliminate the slack Binary variable ($w_{ijk0}$) and its corresponding disaggregated variable ($v_{ijk0}$) from *2E-GDP-HR*. We first combine the last two cardinality rules in *2E-GDP-HR* to obtain (1).

$$w_{ijk0} = 1 - y_{ij} \tag{1}$$

Equating the two variable aggregation constraints in *2E-GDP-HR* and solving for $v_{ijk0}$ gives (2).

$$v_{ijk0} = \sum_{j \in J_i} u_{ij} - \sum_{l \in L_{ijk}} v_{ijkl} \tag{2}$$

Substituting (1) and (2) into the bounding constraint for $v_{ijk0}$ gives (3), which can be rearranged into (4).

$$x^{LB} \cdot (1 - y_{ij}) \leq \sum_{j \in J_i} u_{ij} - \sum_{l \in L_{ijk}} v_{ijkl} \leq x^{UB} \cdot (1 - y_{ij}) \tag{3}$$

$$u_{ij} + \sum_{j' \in J_i : j' \neq j} u_{ij'} - x^{UB} \cdot (1 - y_{ij}) \leq \sum_{l \in L_{ijk}} v_{ijkl} \leq u_{ij} + \sum_{j' \in J_i : j' \neq j} u_{ij'} - x^{LB} \cdot (1 - y_{ij}) \tag{4}$$

Summing the bounding constraint for $x_{ij}$ over $j' \in J_i$ for $j' \neq j$, results in (5). Using the cardinality rule $\sum_{j \in J_i} y_{ij} = 1$, (5) can be written as given in (6), which has two parts, (6a) and (6b). Substituting these into (4) proves that (4) is a relaxation of the disaggregation constraint in *2L-GDP-HR* ($\sum_{l \in L_{ijk}} v_{ijkl} = u_{ij}$, which can be written as $u_{ij} \leq \sum_{l \in L_{ijk}} v_{ijkl} \leq u_{ij}$).

$$\sum_{j' \in J_i : j' \neq j} x^{LB} \cdot y_{ij'} \leq \sum_{j' \in J_i : j' \neq j} u_{ij'} \leq \sum_{j' \in J_i : j' \neq j} x^{UB} \cdot y_{ij'} \tag{5}$$

$$x^{LB} \cdot (1 - y_{ij}) \leq \sum_{j' \in J_i : j' \neq j} u_{ij'} \leq x^{UB} \cdot (1 - y_{ij}) \tag{6}$$

$$\sum_{j' \in J_i : j' \neq j} u_{ij'} - x^{UB} \cdot (1 - y_{ij}) \leq 0 \tag{6a}$$

$$\sum_{j' \in J_i : j' \neq j} u_{ij'} - x^{LB} \cdot (1 - y_{ij}) \geq 0 \tag{6b}$$

It should also be noted that the cardinality rule on the extracted lower-level decisions in *2E-GDP-HR* ($w_{ijk0} + \sum_{l \in L_{ijk}} w_{ijkl} = 1$) is redundant with respect to the other two cardinality rules. This can be shown by noting that $w_{ijk0}$ acts like a slack variable, which allows writing the mentioned cardinality rule as $\sum_{l \in L_{ijk}} w_{ijkl} \leq 1$. This expression is contained in the first two cardinality rules since $y_{ij} \leq 1$ and $\sum_{l \in L_{ijk}} w_{ijkl} = y_{ij}$. Therefore, the Hull reformulation of the Equivalent Single-Level GDP produces constraints with continuous relaxations that are weaker than those resulting from the Hull reformulation of the Nested GDP, proving that *2L-GDP-HR ⊆ 2E-GDP-HR*.

*QED*

*Theorem 2.* Let *rML-GDP-BM* denote the continuous relaxation of the mixed-integer program (MIP) obtained from a Multi-Level Nested GDP via the Big-M reformulation, and let *rME-GDP-BM* denote the continuous relaxation of the MIP obtained from its respective Equivalent Single-Level GDP representation via the Big-M reformulation. The feasible space of the former is contained within the feasible space of the latter, namely, *rML-GDP-BM ⊆ rME-GDP-BM*, if tight values for the *M* parameters are used.

*Proof.* Without loss of generality, the above theorem is proved by establishing that the Big-M reformulation of the 2-Level Nested GDP model (*r2L-GDP-BM*) is contained in the Big-M reformulation of its Equivalent Single-Level GDP representation (*r2E-GDP-BM*), when tight *M* values are used:

*r2L-GDP-BM ⊆ r2E-GDP-BM*

The Big-M reformulation for the nested GDP model is given in *2L-GDP-BM*, where $M_{ij}$ is the Big-M value for the constraints in the $j^{th}$ disjunct in disjunction $i$, $M'_{ijkl}$ is the Big-M value associated with the upper-level decision on the nested constraints, and $m'_{ijkl}$ is the Big-M value associated with the lower-level decision on the nested constraints. The Big-M reformulation for the Equivalent Single-Level GDP is given in *2E-GDP-BM*, where $M_{ij}$ is the same as in *2L-GDP-BM*, and $m_{ijkl}$ is the Big-M value associated with the extracted lower-level decisions.

$$\begin{aligned}
\min z &= f(x) & \text{(2L-GDP-BM)} \\
s.t. \quad & r(x) \leq 0 \\
& g_{ij}(x) \leq M_{ij} \cdot (1 - y_{ij}) & \forall i \in I, j \in J_i \\
& h_{ijkl}(x) \leq m'_{ijkl} \cdot (1 - w_{ijkl}) + M'_{ijkl} \cdot (1 - y_{ij}) & \forall i \in I, j \in J_i, k \in K_{ij}, l \in L_{ijk} \\
& \sum_{j \in J_i} y_{ij} = 1 & \forall i \in I \\
& y_{ij} = \sum_{l \in L_{ijk}} w_{ijkl} & \forall i \in I, j \in J_i, k \in K_{ij}, L_{ijk} \neq \emptyset \\
& Ay + Bw \leq c
\end{aligned}$$

$$x^{LB} \leq x \leq x^{UB}$$
$$x \in \mathbb{R}^n$$
$$y_{ij} \in \{0,1\} \qquad \forall i \in I, j \in J_i$$
$$w_{ijkl} \in \{0,1\} \qquad \forall i \in I, j \in J_i, k \in K_{ij}, l \in L_{ijk}$$

$$\min z = f(x) \qquad \text{(2L-GDP-BM)}$$
$$\text{s.t.} \quad r(x) \leq 0$$
$$g_{ij}(x) \leq M_{ij} \cdot (1 - y_{ij}) \qquad \forall i \in I, j \in J_i$$
$$h_{ijkl}(x) \leq m_{ijkl} \cdot (1 - w_{ijkl}) \qquad \forall i \in I, j \in J_i, k \in K_{ij}, l \in L_{ijk}$$
$$\sum_{j \in J_i} y_{ij} = 1 \qquad \forall i \in I$$
$$y_{ij} = \sum_{l \in L_{ijk}} w_{ijkl} \qquad \forall i \in I, j \in J_i, k \in K_{ij}, L_{ijk} \neq \emptyset$$
$$w_{ijk0} + \sum_{l \in L_{ijk}} w_{ijkl} = 1 \qquad \forall i \in I, j \in J_i, k \in K_{ij}, L_{ijk} \neq \emptyset$$
$$Ay + Bw \leq c$$
$$x^{LB} \leq x \leq x^{UB}$$
$$x \in \mathbb{R}^n$$
$$y_{ij} \in \{0,1\} \qquad \forall i \in I, j \in J_i$$
$$w_{ijkl} \in \{0,1\} \qquad \forall i \in I, j \in J_i, k \in K_{ij}, l \in L_{ijk}$$

Finding the tightest Big-M values requires solving multiple optimization problems to maximize the value of each constraint function over the complete model's feasible region, or over the corresponding feasible region of the disjunction (Grossmann & Trespalacios, 2013). For the proof we calculate tight Big-M values using only the global constraints or upper-level constraints in the case of the nested constraints. The following mathematical optimization problems are solved to obtain tight $M$ values: (7) for $M_{ij}$, (8a) for $m'_{ijkl}$, (8b) for $M'_{ijkl}$, and (9) for $m_{ijkl}$. It should be noted that $m'_{ijkl}$ accounts for the upper-level constraints $g_{ij}(x) \leq 0$, meaning it is localized to the parent disjunct that it belongs to. $M'_{ijkl}$ subtracts $m'_{ijkl}$ from the traditional Big-M value to ensure that when both upper and lower-level decisions are not selected ($y_{ij} = 0$ and $w_{ijkl} = 0$), the resulting Big-M value is equivalent to the global Big-M value for that constraint.

$$M_{ij} = \max\{g_{ij}(x) \mid r(x) \leq 0, x^{LB} \leq x \leq x^{UB}, x \in \mathbb{R}^n\} \tag{7}$$
$$m'_{ijkl} = \max\{h_{ijkl}(x) \mid r(x) \leq 0, g_{ij}(x) \leq 0, x^{LB} \leq x \leq x^{UB}, x \in \mathbb{R}^n\} \tag{8a}$$
$$M'_{ijkl} = \max\{h_{ijkl}(x) \mid r(x) \leq 0, x^{LB} \leq x \leq x^{UB}, x \in \mathbb{R}^n\} - m'_{ijkl} \tag{8b}$$
$$m_{ijkl} = \max\{h_{ijkl}(x) \mid r(x) \leq 0, x^{LB} \leq x \leq x^{UB}, x \in \mathbb{R}^n\} \tag{9}$$

The proof lies in establishing that the feasible space of *2L-GDP-BM* is contained in *2E-GDP-BM*. The difference between these two models is shown in the highlighted constraints above. It was previously shown that the cardinality rule $w_{ijk0} + \sum_{l \in L_{ijk}} w_{ijkl} = 1$ is redundant (see *Theorem 1*). Thus, the proof is given by establishing that the right-hand-sides of the highlighted Big-M constraints satisfy (10), meaning that the Big-M constraint from *2L-GDP-BM* is contained in the Big-M constraint from *2E-GDP-BM*. Substituting (9) in (8b), results in (11). Substituting (11) in (10) and simplifying the resulting expression produces (12). From the cardinality constraint $\sum_{l \in L_{ijk}} w_{ijkl} = y_{ij}$, it is clear that $w_{ijkl} \leq y_{ij}$, meaning that the expressions in parenthesis in (12) can be dropped without changing the sign on the inequality. Thus, $m'_{ijkl} \leq m_{ijkl}$, which is true considering that (9) is a relaxation of (8a). Therefore, *2L-GDP-BM* $\subseteq$ *2E-GDP-BM*.

$$m'_{ijkl} \cdot (1 - w_{ijkl}) + M'_{ijkl} \cdot (1 - y_{ij}) \leq m_{ijkl} \cdot (1 - w_{ijkl}) \tag{10}$$

$$M'_{ijkl} = m_{ijkl} - m'_{ijkl} \tag{11}$$

$$m'_{ijkl} \cdot (y_{ij} - w_{ijkl}) \leq m_{ijkl} \cdot (y_{ij} - w_{ijkl}) \tag{12}$$

QED

### 3.4. Flattening via Basic Steps

The third approach to modeling hierarchical GDP is to flatten the nested disjunctions by applying sufficiently many basic steps within each disjunction until the nested system is transformed into a system with single-level disjunctions. Consider the simple nested disjunction in (13). This disjunction constraint can be flattened by applying two basic steps to introduce $g_1(x) \leq 0$ into the nested disjunctions, resulting in (14), where $Z_1 = Y_1 \wedge W_1$ and $Z_2 = Y_1 \wedge W_2$.

$$\begin{bmatrix} Y_1 \\ g_1(x) \leq 0 \\ \begin{bmatrix} W_1 \\ h_1(x) \leq 0 \end{bmatrix} \vee \begin{bmatrix} W_2 \\ h_2(x) \leq 0 \end{bmatrix} \end{bmatrix} \vee \begin{bmatrix} Y_2 \\ g_2(x) \leq 0 \end{bmatrix} \tag{13}$$

$$\begin{bmatrix} Z_1 \\ g_1(x) \leq 0 \\ h_1(x) \leq 0 \end{bmatrix} \vee \begin{bmatrix} Z_2 \\ g_1(x) \leq 0 \\ h_2(x) \leq 0 \end{bmatrix} \vee \begin{bmatrix} Y_2 \\ g_2(x) \leq 0 \end{bmatrix} \tag{14}$$

For disjunctions with a single nested disjunction, applying a basic step is quite cheap. However, once there is more than one nested disjunction inside a single disjunct, the number of basic steps required to flatten the hierarchical GDP grows exponentially. Consider the disjunction with two nested disjunctions in (15). Flattening the disjunction is a set covering problem and requires eight basic steps (four for each combination of two disjuncts and four more to introduce $g_1(x) \leq 0$ in the resulting disjuncts) to obtain the equivalent disjunction in (16).

$$\begin{bmatrix} Y_1 \\ g_1(x) \leq 0 \\ \begin{bmatrix} W_1 \\ h_1(x) \leq 0 \end{bmatrix} \vee \begin{bmatrix} W_2 \\ h_2(x) \leq 0 \end{bmatrix} \\ \begin{bmatrix} W_3 \\ h_3(x) \leq 0 \end{bmatrix} \vee \begin{bmatrix} W_4 \\ h_4(x) \leq 0 \end{bmatrix} \end{bmatrix} \vee \begin{bmatrix} Y_2 \\ g_2(x) \leq 0 \end{bmatrix} \tag{15}$$

$$\begin{bmatrix} Y_1 \wedge W_1 \wedge W_3 \\ g_1(x) \leq 0 \\ h_1(x) \leq 0 \\ h_3(x) \leq 0 \end{bmatrix} \vee \begin{bmatrix} Y_1 \wedge W_1 \wedge W_4 \\ g_1(x) \leq 0 \\ h_1(x) \leq 0 \\ h_4(x) \leq 0 \end{bmatrix} \vee \begin{bmatrix} Y_1 \wedge W_2 \wedge W_3 \\ g_1(x) \leq 0 \\ h_2(x) \leq 0 \\ h_3(x) \leq 0 \end{bmatrix} \vee \begin{bmatrix} Y_1 \wedge W_2 \wedge W_4 \\ g_1(x) \leq 0 \\ h_2(x) \leq 0 \\ h_4(x) \leq 0 \end{bmatrix} \vee \begin{bmatrix} Y_2 \\ g_2(x) \leq 0 \end{bmatrix} \quad (16)$$

Generalizing this to the notation of *2L-GDP*, a disjunction with $k \in K_{ij}$ nested disjunctions, each of which has $l \in L_{ijk}$ disjuncts, requires the number of basic steps given in (17), where the notation $\binom{a}{1}$ is the binomial coefficient (choose 1 from a group with $a$ elements). The coefficient 2 accounts for introducing $g_{ij}(x) \leq 0$ into each of the resulting disjuncts, and can be replaced by the $1 + |g_{ij}(x)|$ if $g_{ij}(x)$ represents a vector of functions, where $|g_{ij}(x)|$ is the number of functions within $g_{ij}(x)$. A hybrid approach is also possible, where some basic steps are performed and then the resulting nested disjunction is flattened as in the Equivalent Single-Level GDP approach. However, as the number of nested disjunctions increases, this hybrid approach results in many more disjunctions than those given in (17). Although flattening via basic steps may produce models that are tighter than the inside-out reformulation of the Nested GDP, the combinatorial growth of such systems makes this approach prohibitive for multi-level decision systems with multiple disjuncts in each nested disjunction.

$$2 \cdot \binom{|L_{ij1}|}{1} \cdot \binom{|L_{ij2}|}{1} \cdot \ldots \cdot \binom{|L_{ij|K_{ij}|}|}{1} = 2 \cdot \prod_{k \in K_{ij}} |L_{ijk}| \quad (17)$$

## 4. Examples

Each of the examples in this section are implemented in the Julia programming language (version 1.8.3) (Bezanson et al., 2017) using various packages within the ecosystem. These include *JuMP* (version 1.6.0) (Dunning et al., 2017) for modeling mathematical programs, *DisjunctiveProgramming* (version 0.3.6) (Perez et al., 2023) for reformulating GDPs (both nested and single-level) into MIPs, and *Polyhedra* (version 0.7.5) (Legat et al., 2021) for projecting mathematical programming models onto 2D space (see **Example 4.1**). For the numerical examples (**Examples 4.2** and **4.3**), the reformulated MI(N)LP models are solved on a Windows PC with 16 GB of RAM and an Intel-i7 3.60 GHz processor. CPLEX (version 22.1.0) is used as the MILP solver and BARON (version 23.1.5) as the MINLP solver.

### 4.1. Graphical Example of Model Tightness

Consider the Nested GDP constraint system given in (18), which can be expressed as the Equivalent Single-Level GDP in (19), or the Flattened GDP (via basic steps) in (20). Each of these models is reformulated into a MIP using the Big-M reformulation, with both a loose (large) *M* value and a tight *M* value, and the Hull reformulation. Their continuous relaxations are then projected onto the $x_1, x_2$ plane in **Figure 4.1**. The projections show that flattening via basic steps is advantageous when the Hull reformulation is performed, but not necessarily when the Big-M reformulation is performed with a tight *M* value. In the latter case, the lack of information regarding the system hierarchy results in a Big-M reformulation that is equivalent to taking the Big-M region between $W_1$, $W_2$, and $Y_2$, which is worse than taking the Big-M region of $W_1$ and $W_2$ and intersecting it with the Big-M region of $Y_1$ and $Y_2$.

Explicitly preserving the hierarchical relationship in the Nested GDP representation reduces the feasible region of the continuous relaxation more than when the Equivalent Single-Level GDP representation is used. This is observed in both the Tight-M (Big-M reformulation with a tight *M*) and Hull reformulation cases. Furthermore, in this example the Tight-M reformulation of the Nested GDP model produces the same relaxation as the Hull reformulation of the Equivalent Single-Level GDP model with only a fraction of the model size (see **Table 4.1**). It should also be noted that, the convex hull of the system is obtained when either the Hull reformulation is applied to the Nested GDP or when it is applied to the Flattened GDP. As a result, the continuous relaxation of either formulation will yield the optimum.

$$\begin{bmatrix} Y_1 \\ 1 \leq x_1 \leq 3 \\ 4 \leq x_2 \leq 6 \\ \begin{bmatrix} W_1 \\ 1 \leq x_1 \leq 2 \\ 5 \leq x_2 \leq 6 \end{bmatrix} \vee \begin{bmatrix} W_2 \\ 2 \leq x_1 \leq 3 \\ 4 \leq x_2 \leq 5 \end{bmatrix} \end{bmatrix} \vee \begin{bmatrix} Y_2 \\ 8 \leq x_1 \leq 9 \\ 1 \leq x_2 \leq 2 \end{bmatrix} \quad (18)$$

$$\Xi(1, \{Y_1, Y_2\})$$
$$\Xi(Y_1, \{W_1, W_2\})$$

$$\begin{bmatrix} Y_1 \\ 1 \leq x_1 \leq 3 \\ 4 \leq x_2 \leq 6 \end{bmatrix} \vee \begin{bmatrix} Y_2 \\ 8 \leq x_1 \leq 9 \\ 1 \leq x_2 \leq 2 \end{bmatrix}$$

$$\begin{bmatrix} W_1 \\ 1 \leq x_1 \leq 2 \\ 5 \leq x_2 \leq 6 \end{bmatrix} \vee \begin{bmatrix} W_2 \\ 2 \leq x_1 \leq 3 \\ 4 \leq x_2 \leq 5 \end{bmatrix} \vee \begin{bmatrix} W_3 \\ 1 \leq x_1 \leq 9 \\ 1 \leq x_2 \leq 6 \end{bmatrix} \quad (19)$$

$$\Xi(1, \{Y_1, Y_2\})$$
$$\Xi(1, \{W_1, W_2, W_3\})$$
$$\Xi(Y_1, \{W_1, W_2\})$$

$$\begin{bmatrix} YW_{11} \\ 1 \leq x_1 \leq 3 \\ 4 \leq x_2 \leq 6 \\ 1 \leq x_1 \leq 2 \\ 5 \leq x_2 \leq 6 \end{bmatrix} \vee \begin{bmatrix} YW_{12} \\ 1 \leq x_1 \leq 3 \\ 4 \leq x_2 \leq 6 \\ 2 \leq x_1 \leq 3 \\ 4 \leq x_2 \leq 5 \end{bmatrix} \vee \begin{bmatrix} Y_2 \\ 8 \leq x_1 \leq 9 \\ 1 \leq x_2 \leq 2 \end{bmatrix} \quad (20)$$

$$\Xi(1, \{YW_{11}, YW_{12}, Y_2\})$$

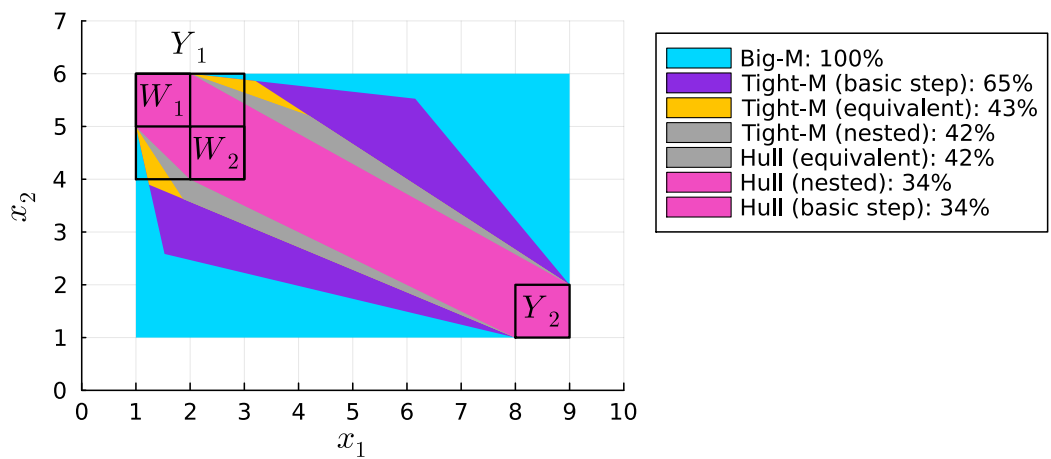

**Figure 4.1.** Projections of the continuous relaxations of Example 4.1 for different reformulations (Big-M = Big-M Reformulation, Tight-M = Big-M Reformulation with tight *M* values, Hull = Hull Reformulation; basic step = Flatten GDP via basic steps, equivalent = Equivalent Single-Level GDP, nested = 2-Level Nested GDP). Projection areas, relative to Big-M are indicated in %.

**Table 4.1.** Model sizes and projection areas for Example 4.1

| Approach | Binary Variables | Continuous Variables | Constraints | Feasible Area | Relative Area |
|---|---|---|---|---|---|
| Big-M | 5 | 2 | 28 | 40.0 | 100% |
| Tight-M (basic step) | 3 | 2 | 28 | 26.0 | 65% |
| Tight-M (equivalent) | 5 | 2 | 28 | 17.3 | 43% |
| Tight-M (nested) | 4 | 2 | 26 | 16.7 | 42% |
| Hull (equivalent) | 5 | 12 | 72 | 16.7 | 42% |
| Hull (nested) | 4 | 10 | 64 | 13.5 | 34% |
| Hull (basic step) | 3 | 8 | 54 | 13.5 | 34% |

## 4.2. Process Superstructure Optimization with Technology Selection and Scheduling Decisions

Consider the superstructure optimization problem for a plant that is to produce and sell material *D* (see **Figure 4.2.1**). Material *D* can be produced from material *C* (reaction: *C → D*), which can be purchased from a third party or produced from material *B* (reaction: *B → C*), which can in turn be purchased or produced from material *A* (reaction: *A → B*). The plant has two types of multipurpose reactors, each with a backup unit, that can be used for the material transformation steps (see **Figure 4.2.2**). Each of these has a maximum installed capacity of 100 kg. Up to one tank for each material in the system can be installed for storage with a maximum installed capacity of 300 kg. There are two candidate chemical processes to perform each material transformation step, giving a total of six processes in the process superstructure. There are two potential technologies (catalysts) that can be used in each process, each with a unique cost and yield, giving a total of 12 candidate process-catalysts combinations in the system. The plant process and equipment superstructures are given in **Figures 4.2.1** and **4.2.2**, respectively. The former illustrates the candidate processes in the superstructure in the state-task network representation (Kondili et al., 1993). The latter depicts the equipment options (reactor type and units, and tanks) in the superstructure.

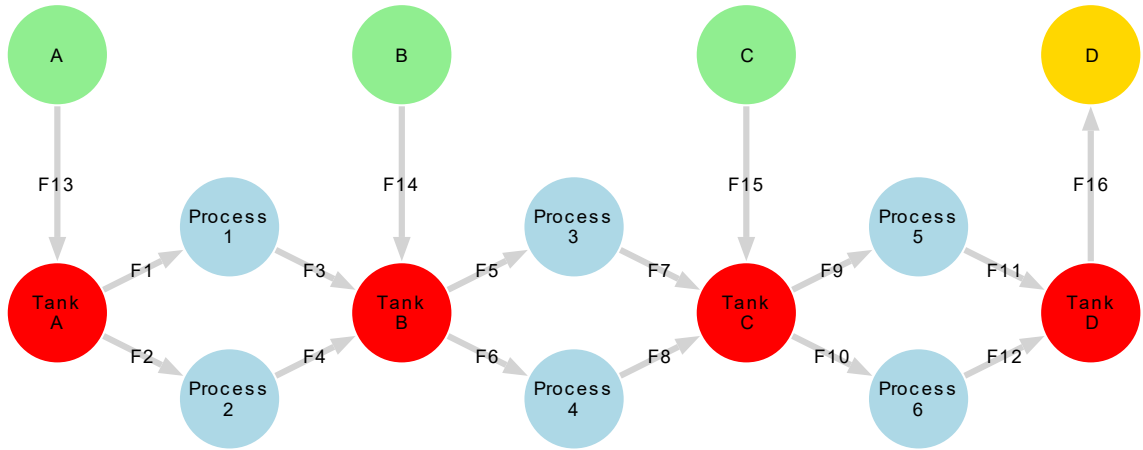

**Figure 4.2.1.** Process superstructure for Example 4.2 with 4 materials, 6 processes, 4 tanks, and 16 streams.

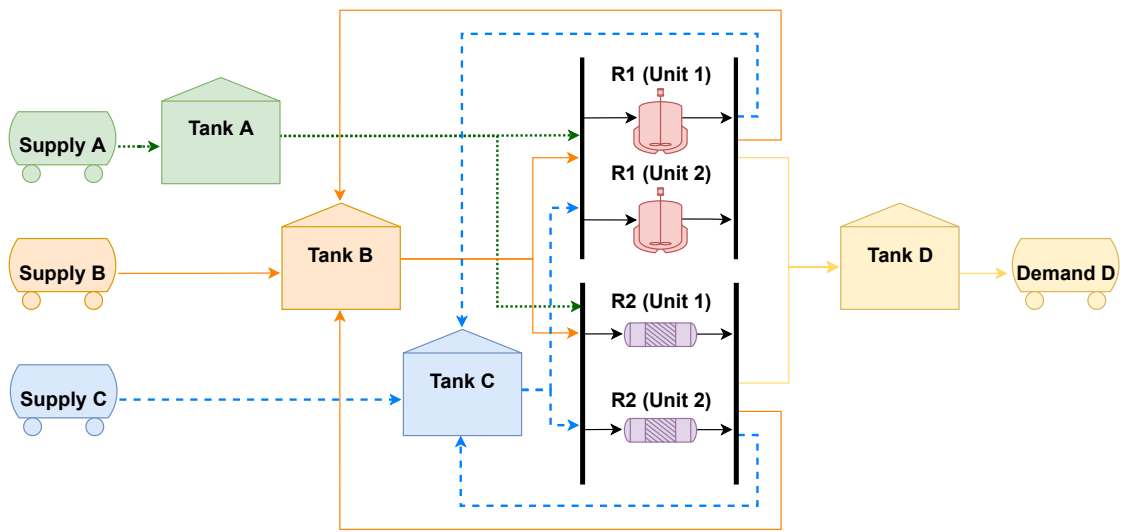

**Figure 4.2.2.** Equipment superstructure (process flow diagram) for Example 4.2 with 4 tanks and 2 reactor types, each with 2 identical units.

The objective of the optimization problem is to maximize system profit over a 30-day schedule by making the following decisions:

- Which material storage tanks to install.
- How many shared reactors to install.
- Which processes to install for each material transformation step.
    - Which technologies (catalysts) to use in each of the selected processes.
    - Which reactor type to use in each of the selected processes.
        - How many reactors to operate in each time period.
        - How much to produce in each batch of material.
- How much material to purchase for *A*, *B*, and *C* in each time period.

The hierarchy of these decisions is indicated by the bullet indentation above. Thus, the technology and reactor type selections are second-level decisions, and the operating schedule and batch sizes are third-level decisions. For simplicity, any changeover or setup times are not considered.

*Model:* The model for this system consists of the following linear constraints. Resource balances are enforced around each resource $k$ at timepoint $t$ with the global constraints in (21) and (22). The level of material at each tank, $L_{k,t}$, is updated based on the material flowing in and out of the tank (material balance). The availability of each reactor, $R_{k,t}$, is updated based on the reactor usage, $\Delta R_{i,k,t}$. A reactor unit is locked (unavailable) when it begins a processing task $i$ at time $t$. At time $t + \tau_i$, the processing task ends ($\tau_i$ is the duration), and the reactor unit is released (becomes available). The values used for the task durations, $\tau_i$, are $\tau_i = 5 \ \forall i \in \{1,4,5,6\}$, $\tau_2 = 3$, and $\tau_3 = 4$ (days). For greater detail on resource balances, the reader is referenced to the review paper on the resource-task network by Perez et al. (2022).

$$\overbrace{L_{k,t}}^{tank\ level} = L_{k,t-1} + \overbrace{\sum_{s \in S_k^{in}} F_{s,t}}^{inflow} - \overbrace{\sum_{s \in S_k^{out}} F_{s,t}}^{outflow} \qquad \forall k \in K^{tank}, t \in T \qquad (21)$$

$$\overbrace{R_{k,t}}^{reactor\ availability} = R_{k,t-1} + \sum_{i \in I} \left( \overbrace{\Delta R_{i,k,t-\tau_i}}^{reactors\ released} - \overbrace{\Delta R_{i,k,t}}^{reactors\ locked} \right) \qquad \forall k \in K^{react}, t \in T \qquad (22)$$

The decision to install a resource (tank or reactor) is governed by the disjunctions in (23) and (24), where the decision is to determine how many units $u$ to install. In this example, $U_k = \{0,1,2\}$ for each reactor type (at most 2 identical units can be installed for each reactor type $k$), and $U_k = \{0,1\}$ for each tank (at most 1 tank can be installed for each material). The installation cost, $CI_k$, is calculated as the sum of a fixed charge, $\alpha_k$, and a variable cost coefficient, $\beta_k$, times the total resource capacity. If no units are installed ($u = 0$), the installation cost and resource capacity, $Q_k$, drop to zero. (23) and (24) also set the initial condition for the resource availability, $L_{k,0}$ and $R_{k,0}$: if installed, tanks are full, and all reactor units are available, respectively. (23) also tracks the slack on the tank level at the final timepoint $|T|$, $\hat{L}_k$, which is penalized in the objective function to reduce the likelihood of depleting the inventory at the end of the scheduling horizon (see (47)). These constraints ensure that the schedule obtained is a feasible schedule for normal operation with monthly cycles. For startup operations the optimal schedule can be obtained by fixing the design decisions and rerunning the model with the initial tank levels set to zero. The cardinality constraint (25) ensures that exactly one of the disjuncts is selected. The values for the cost coefficients are given in **Table 4.2.1**. Since the plant lifetime is greater than the scheduling horizon, resource installation costs coefficients have been scaled to the appropriate order of magnitude. Installation costs for pipelines between tanks and reactors are assumed to be negligible.

$$\bigvee_{u \in U_k \setminus \{0\}} \begin{bmatrix} X_{k,u} \\ CI_k = \alpha_k + \beta_k \cdot \overbrace{u \cdot Q_k}^{installed\ capacity} \\ L_{k,0} = u \cdot Q_k \\ L_{k,|T|} + \hat{L}_k = u \cdot Q_k \end{bmatrix} \bigvee \begin{bmatrix} X_{k,0} \\ CI_k = 0 \\ Q_k = 0 \\ L_{k,t} = 0 \ \forall t \in \{0\} \cup T \\ \hat{L}_k = 0 \end{bmatrix} \qquad \forall k \in K^{tank} \qquad (23)$$

$$\bigvee_{u \in U_k \setminus \{0\}} \begin{bmatrix} X_{k,u} \\ CI_k = \alpha_k + \beta_k \cdot \overbrace{u \cdot Q_k}^{\text{installed capacity}} \\ R_{k,0} = u \end{bmatrix} \vee \begin{bmatrix} X_{k,0} \\ CI_k = 0 \\ Q_k = 0 \\ R_{k,0} = 0 \;\; \forall t \in \{0\} \cup T \end{bmatrix} \qquad \forall k \in K^{react} \quad (24)$$

$$\Xi(1, X_{k,u} \;\; \forall u \in U_k) \qquad \forall k \in K \quad (25)$$

**Table 4.2.1.** Fixed and variable cost coefficients for the installation cost of each resource.

| Resource $(k)$ | Fixed Cost Coefficient $(\alpha_k)$ | Variable Cost Coefficient $(\beta_k)$ |
|---|---|---|
| Tank Material A | $0.406 | $0.011/kg |
| Tank Material B | $0.069 | $0.070/kg |
| Tank Material C | $0.862 | $0.029/kg |
| Tank Material D | $0.086 | $0.003/kg |
| Reactor Type 1 | $0.662 | $0.054/kg |
| Reactor Type 2 | $0.116 | $0.090/kg |

The multi-level disjunction in (26) represents the decision to install process $i$ or not. When installed, the total batch size, $B_{i,t}$, is equal to the flow entering the process at time $t$. There are two nested disjunctions if a process is installed. The first of these relates to which reactor type $k$ is assigned to the process, $W_{i,k}$. The second one pertains to which technology (catalyst) is used for that particular process, $\widehat{W}_{i,j}$. Once a reactor type is assigned, the per unit batch size, $\widehat{B}_{i,t}$, is bounded by the installed capacity of each unit, $Q_k$, and the operating cost, $CO_{i,t}$, is proportional to the total batch size with a cost coefficient $\gamma_{i,k}$ (given in **Table 4.2.2**). The nested technology selection disjunction specifies the amount of material leaving the process when the batch is completed. This is governed by the yield, $\nu$, which is specific to the technology $j$ (given in **Table 4.2.3**). There is then a third-level set of disjunctions inside the reactor type assignment disjunction, which determines the number of units, $u$, that are used for a batch at time, $t$, $N_{i,k,t,u}$. The number of units selected indicates the number of units that are locked at time $t$ and is also used to determine the total batch size from the per unit batch size. Note that for this system, it is assumed that if multiple units are used, their loads are equally distributed. Finally, when a process is not installed ($\neg Y_i$), all pertinent variables are set to zero, and the reactor capacity is only bounded by the maximum allowed capacity. The cardinality rules in (27)-(29) are the linking constraints between the different levels of this multi-level disjunction.

$$\begin{bmatrix} Y_i \\ F_{s,t} = B_{i,t} \ \forall s \in S_i^{in} \\ \bigvee_{k \in K^{react}} \begin{bmatrix} W_{i,k} \\ \hat{B}_{i,t} \leq Q_k \ \forall t \in T \\ CO_{i,t} = \gamma_{i,k} \cdot B_{i,t} \ \forall t \in T \\ \bigvee_{u \in U_k} \begin{bmatrix} N_{i,k,t,u} \\ \Delta R_{i,k,t} = u \\ B_{i,t} = \underbrace{u \cdot \hat{B}_{i,t}}_{\substack{total\ batch \\ size}} \end{bmatrix} \forall t \in T \end{bmatrix} \\ \bigvee_{j \in J} \begin{bmatrix} \widehat{W}_{i,j} \\ F_{s,t+\tau_i} = \underbrace{\nu_{i,j} \cdot B_{i,t}}_{\substack{batch \\ yield}} \ \forall s \in S_i^{out}, t \in T \end{bmatrix} \end{bmatrix} \bigvee \begin{bmatrix} \neg Y_i \\ F_{s,t} = 0 \ \forall s \in S_i^{in} \cup S_i^{out} \\ \hat{B}_{i,t}, B_{i,t}, CO_{i,t} = 0 \ \forall t \in T \\ \Delta R_{i,k,t} = 0 \ \forall k \in K^{react}, t \in T \\ Q_k \leq Q_k^{UB} \ \forall k \in K^{react} \end{bmatrix} \quad \forall i \in I \quad (26)$$

$$\Xi(Y_i, W_{i,k} \ \forall k \in K^{react}) \qquad \forall i \in I \quad (27)$$

$$\Xi(W_{i,k}, N_{i,k,t,u} \ \forall u \in U_k) \qquad \forall i \in I, k \in K^{react}, t \in T \quad (28)$$

$$\Xi(Y_i, \widehat{W}_{i,j} \ \forall j \in J) \qquad \forall i \in I \quad (29)$$

**Table 4.2.2.** Operating cost parameter, $\gamma_{i,k}$ ($/kg), for each process $i$ and reactor type $k$ combination.

| Reactor Type ($k$) | Process ($i$) | | | | | |
|---|---|---|---|---|---|---|
| | 1 | 2 | 3 | 4 | 5 | 6 |
| 1 | 0.258 | 0.339 | 0.425 | 0.905 | 0.745 | 0.156 |
| 2 | 0.575 | 0.454 | 0.017 | 0.496 | 0.636 | 0.087 |

**Table 4.2.3.** Production yield parameter, $\nu_{i,j}$, for each process $i$ and technology $j$ combination.

| Technology ($j$) | Process ($i$) | | | | | |
|---|---|---|---|---|---|---|
| | 1 | 2 | 3 | 4 | 5 | 6 |
| 1 | 42.6% | 57.4% | 44.3% | 62.1% | 86.3% | 51.8% |
| 2 | 76.7% | 13.4% | 8.7% | 35.1% | 19.3% | 11.7% |

An additional logic proposition must be included to ensure that if a process $i$ is triggered on reactor type $k$ at time $t$ with $u$ units ($N_{i,k,t,u} = True$), the reactor type $k$ must have been installed with at least $u$ units ($\exists u' \in U_k: u' \geq u, X_{k,u'} = True$). For example, if $N_{i,k,t,1} = True$, then either $X_{k,1} = True$ or $X_{k,2} = True$ (one or two units must have been installed when the plant was built). This condition is enforced with (30).

$$N_{i,k,t,u} \Rightarrow \bigvee_{u' \in U_k: u' \geq u} X_{k,u'} \qquad \forall i \in I, k \in K^{react}, t \in T, u \in U_k \setminus \{0\} \quad (30)$$

The variable bounds and domains are given in (31)-(33) and (35)-(46). The upper bound resource capacities are, $Q_k^{UB} = 300kg \ \forall k \in K^{tank}$ and $Q_k^{UB} = 100kg \ \forall k \in K^{react}$. The initialization constraint in (34) is used to ensure that there is no flow leaving a reactor in the first $\tau_i$ periods since it is assumed that all reactors are idle at the beginning of the scheduling horizon.

$$0 \leq B_{i,t} \leq \sum_{k \in K^{react}} (|U_k| - 1) \cdot Q_k^{UB} \qquad \forall i \in I, t \in T \qquad (31)$$

$$0 \leq \hat{B}_{i,t} \leq \max(Q_k^{UB} \; \forall k \in K^{react}) \qquad \forall i \in I, t \in T \qquad (32)$$

$$0 \leq F_{s,t} \leq F_s^{UB} \qquad \forall s \in S, t \in T \qquad (33)$$

$$F_{s,t} = 0 \qquad \forall i \in I, s \in S_i^{out}, t \in \{1, \ldots, \tau_i\} \qquad (34)$$

$$0 \leq CI_k \leq \alpha_k + \beta_k \cdot (|U_k| - 1) \cdot Q_k^{UB} \qquad \forall k \in K \qquad (35)$$

$$0 \leq L_{k,t} \leq (|U_k| - 1) \cdot Q_k^{UB} \qquad \forall k \in K^{tank}, t \in T \qquad (36)$$

$$0 \leq \hat{L}_k \leq (|U_k| - 1) \cdot Q_k^{UB} \qquad \forall k \in K^{tank} \qquad (37)$$

$$0 \leq CO_{i,t} \leq \sum_{k \in K^{react}} \gamma_{i,k} \cdot (|U_k| - 1) \cdot Q_k^{UB} \qquad \forall i \in I, t \in T \qquad (38)$$

$$0 \leq Q_k \leq Q_k^{UB} \qquad \forall k \in K^{react} \qquad (39)$$

$$0 \leq R_{k,t} \leq |U_k| - 1 \qquad \forall k \in K^{react}, t \in T \qquad (40)$$

$$0 \leq \Delta R_{i,k,t} \leq |U_k| - 1 \qquad \forall i \in I, k \in K^{react}, t \in T \qquad (41)$$

$$N_{i,k,t,u} \in \{True, False\} \qquad \forall i \in I, k \in K^{react}, t \in T, u \in U_k \qquad (42)$$

$$W_{i,k} \in \{True, False\} \qquad \forall i \in I, k \in K^{react} \qquad (43)$$

$$\widehat{W}_{i,j} \in \{True, False\} \qquad \forall i \in I, j \in J \qquad (44)$$

$$X_{k,u} \in \{True, False\} \qquad \forall k \in K, u \in U_k \qquad (45)$$

$$Y_i \in \{True, False\} \qquad \forall i \in I \qquad (46)$$

The objective of this optimization problem is to maximize profit, as given by (47), where $p_s$ is the price of each material ($p_A = -\$1/kg$, $p_B = -\$7/kg$, $p_C = -\$8/kg$, and $p_D = \$10/kg$). The tank level slacks are penalized with a penalty coefficient equal to the absolute value of the material price.

$$\max \sum_{t \in T} \left( \overbrace{\sum_{s \in S^{mat}} p_s \cdot F_{s,t}}^{\text{material sales/purchases}} - \overbrace{\sum_{i \in I} CO_{i,t}}^{\text{operating costs}} \right) - \overbrace{\sum_{k \in K} CI_k}^{\text{installation costs}} - \overbrace{\sum_{k \in K^{tank}} |p_k| \cdot \hat{F}_k}^{\text{tank slack penalties}} \qquad (47)$$

The resulting model is the linear Nested GDP given in (21)-(47). This hierarchical model is reformulated into a mixed-integer linear program (MILP) using both Big-M (with both loose and tight *M* values) and Hull reformulations. The hierarchical GDP model is also transformed into its Equivalent Single-Level GDP and reformulated with both Big-M and Hull methods. The flattening via basic steps is not performed because the scheduling decisions are made up of $|T|$ nested disjunctions in the third level of nesting, which makes the flattening procedure impractical due to the exponential growth in disjunctions.

The optimum solution yields a cumulative profit of $2,085. The process network and equipment network designs are given in **Figures 4.2.3** and **4.2.4**, respectively. The Gantt charts for procurement/sales and production are shown in **Figures 4.2.5** and **4.2.6**, respectively. The tank levels are displayed in **Figure 4.2.7**. The optimal design requires the installation of Processes 1, 3-5; Tanks B and C; and both reactor types, each with two units available. Reactors of type 1 focus almost exclusively on Process 1 with Technology 2, with one batch of Process 3 (Technology 1). Rectors of type 2 are used for Processes 4 and 5, each using Technology 1. Procurement of A occurs every 5 days, with sales of D typically spaced out every 10 days. By the end of the scheduling horizon, both tank levels have been restored to their initial levels (full).

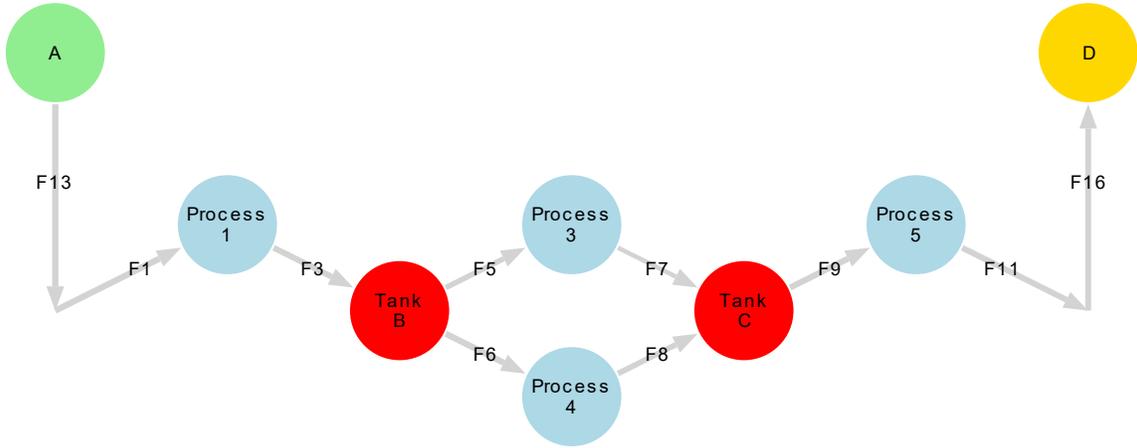

**Figure 4.2.3.** Optimal process network design (edge thickness is proportional to the maximum flow on that line).

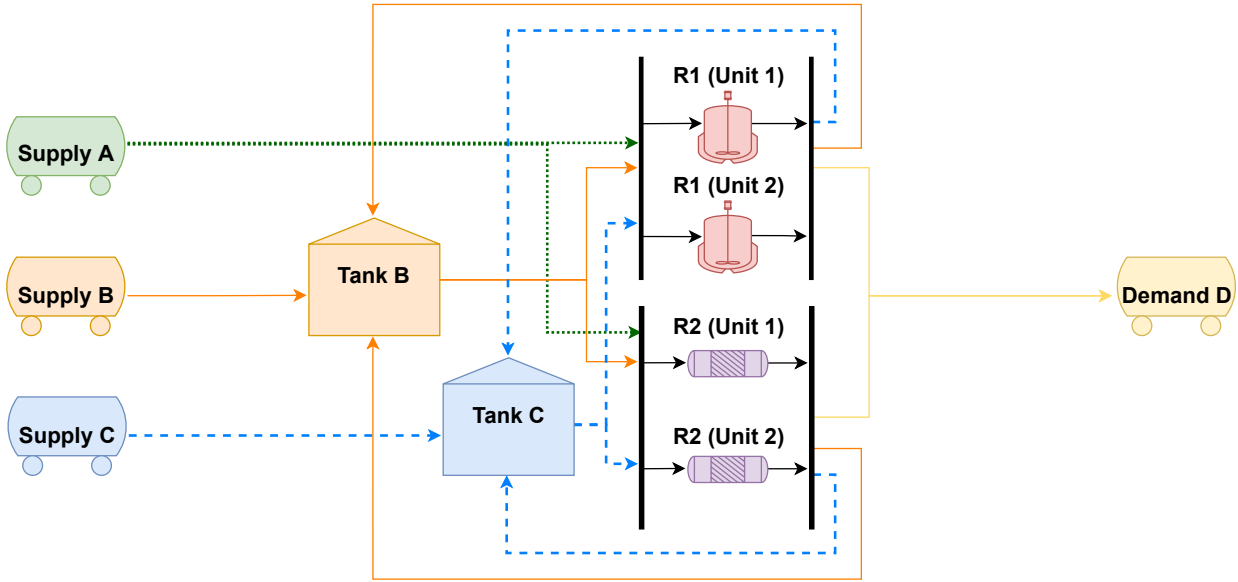

**Figure 4.2.4.** Optimal equipment network design.

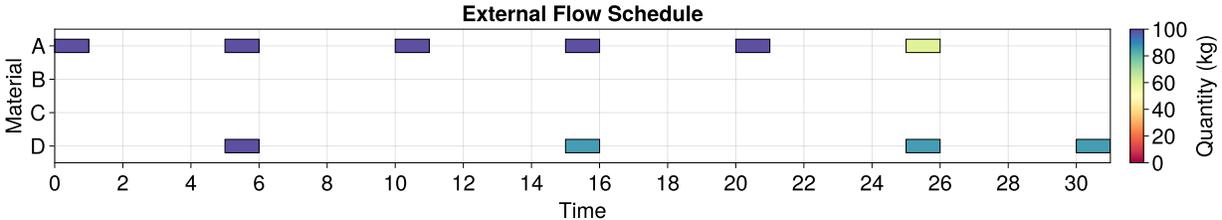

**Figure 4.2.5.** Material procurement and sales schedule.

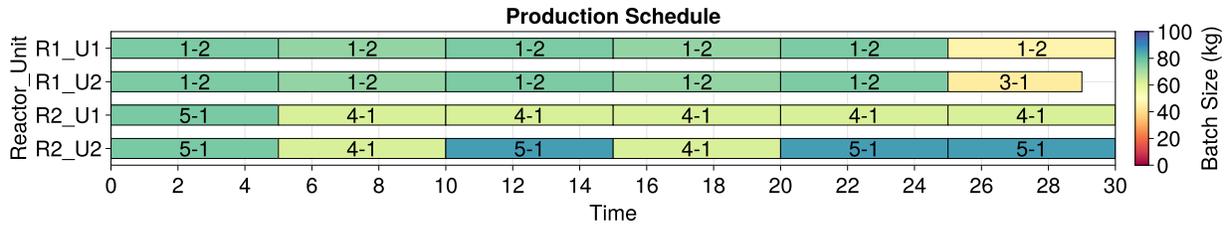

**Figure 4.2.6.** Plant operations schedule (text in each bar, *i-j*, indicates process number *i* and technology number *j* for that event).

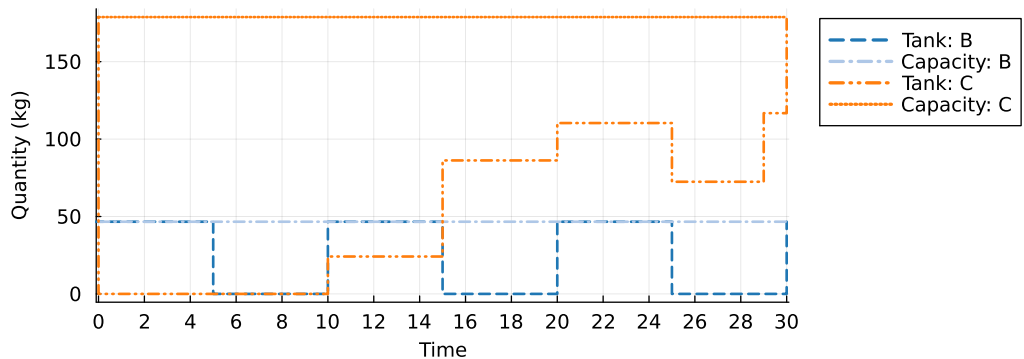

**Figure 4.2.7.** Amount of material in each tank (and maximum tank level) throughout the scheduling horizon.

The model sizes and computational statistics for each of the reformulated MILP models is given in **Table 4.2.4**. Each model is solved three times: 1) relaxing the integrality constraints (LP relaxation), 2) solving the MILP with presolve and heuristics disabled in CPLEX to compare the effects of the model formulation size and tightness on the computational performance, and 3) solving the MILP with presolve and heuristics enabled. As can been observed, all formulations, except the Hull reformulation of the Nested GDP, have a poor continuous relaxation with very large relaxation gaps. The Hull reformulated Nested GDP, on the other hand, has a tight relaxation with a 9% relative gap. When presolve and heuristics are disabled in CPLEX, both the Big-M and Tight-M reformulations cannot solve the problem to optimality within the allotted time limit of 3,600 seconds. CPLEX is able to reduce the optimality gap more in the Tight-M models, than in their Big-M counterparts, with a greater gap reduction when the Nested versions are used. The Tight-M models result in solutions that are approximately 10% lower than the optimum. It is also observed that the Equivalent Single-Level models result in better feasible solutions despite the larger optimality gaps. The poor solution found with the Nested Big-M model is likely due in part to the fact that the nested constraints end up having two very large *M* values when reformulated, making them less tight than their Single-Level counterparts. Interestingly, when presolve and heuristics are enabled in CPLEX, it is able to solve the models with large *M* values to optimality, but not those with tighter *M* values.

The MILP model obtained by applying the Hull reformulation to the Nested GDP model outperforms the other models, finding the optimum in approximately half of the time required relative to its Equivalent Single-Level counterpart, requiring less nodes to be explored and less cuts to be applied. This superior performance is likely due to the reduced model size and the LP relaxation tightness. The Hull reformulated

Nested GDP results in a model that has fewer binary variables (25% and 4% less, before and after presolve, respectively), continuous variables (10% and 26% less, before and after presolve, respectively), and constraints (5% and 25% less, before and after presolve, respectively) than its Equivalent Single-Level counterpart.

**Table 4.2.4.** Model sizes and computational statistics of the MILP models resulting from the Big-M reformulations (using both loose and tight $M$ values) and Hull reformulations of the GDP models in Example 4.2.

|  | Big-M Reformulation | | Tight-M Reformulation | | Hull Reformulation | |
| --- | --- | --- | --- | --- | --- | --- |
|  | *Equivalent* | *Nested* | *Equivalent* | *Nested* | *Equivalent* | *Nested* |
| **LP Relaxation** | | | | | | |
| *Relaxation Solution* | $81,000 | $81,000 | $80,999 | $80,999 | $68,044 | $2,268 |
| *Relaxation Gap* | 3785% | 3785% | 3785% | 3785% | 3163% | 9% |
| **MIP Solution**[a] | | | | | | |
| *MIP Solution* | $1,646 | $784 | $1,874 | $1,838 | $2,085 | $2,085 |
| *MIP Gap* | 3888% | 1454% | 303% | 177% | 0% | 0% |
| *Nodes Explored* | 721,970 | 1,470,170 | 970,567 | 1,564,375 | 80,388 | 45,094 |
| *Cuts Applied* | 2,641 | 1,080 | 1,958 | 1,253 | 852 | 94 |
| *CPU Time (s)* | 3,602 | 3,603 | 3,603 | 3,604 | 878 | 411 |
| **MIP Solution** | | | | | | |
| *MIP Solution* | $2,085 | $2,085 | $2,085 | $2,083 | $2,085 | $2,085 |
| *MIP Gap* | 0% | 0% | 42% | 29% | 0% | 0% |
| *Nodes Explored* | 129,716 | 593,418 | 284,619 | 462,714 | 35,920 | 31,993 |
| *Cuts Applied* | 410 | 381 | 629 | 762 | 746 | 22 |
| *CPU Time (s)* | 427 | 1,031 | 3,612 | 3,603 | 81 | 38 |
| **Original Model Size** | | | | | | |
| *Binary Variables* | 1,550 | 1,166 | 1,550 | 1,166 | 1,550 | 1,166 |
| *Continuous Variables* | 1,634 | 1,634 | 1,634 | 1,634 | 11,850 | 10,614 |
| *Constraints* | 14,817 | 14,049 | 14,817 | 14,049 | 59,195 | 56,243 |
| **Reduced Model Size**[b] | | | | | | |
| *Binary Variables* | 1,158 | 1,156 | 1,162 | 1,156 | 1,158 | 1,115 |
| *Continuous Variables* | 1,409 | 1,409 | 1,409 | 1,409 | 4,468 | 3,302 |
| *Constraints* | 7,443 | 6,969 | 7,849 | 7,843 | 8,166 | 6,163 |

[a]No presolve; no heuristics
[b]After presolve (CPLEX applies presolve 3 times for Big-M and 2 times for others)

## 4.3. Multi-period Design and Expansion Planning

Example 4.3 is based on Example 1 in the work by van den Heever and Grossmann (1999), which consists of an integrated superstructure optimization problem with long term operational and expansion planning. The problem has three potential processes (1, 2, and 3), each with its dedicated processing unit, and three materials (*A*, *B*, and *C*) as shown in **Figure 4.3.1**. Material *C* is the final product (price: $10,800/ton) and is produced from Material *B* in Process 1. Material *B* can be purchased externally (cost: $7,000/ton) or produced from Material A (cost: $1,800/ton) in either Process 2 or Process 3. It is assumed that each process includes any required separation steps, such that the respective exit streams are single-component streams containing the pure product of each process. The objective here is to minimize cost (maximize profit) by making the following decisions:

- Which processes should be used.
- Which processes to operate in each period.
- Which processes to undergo a capacity expansion in each period.
- How much new processing capacity to install in each period.

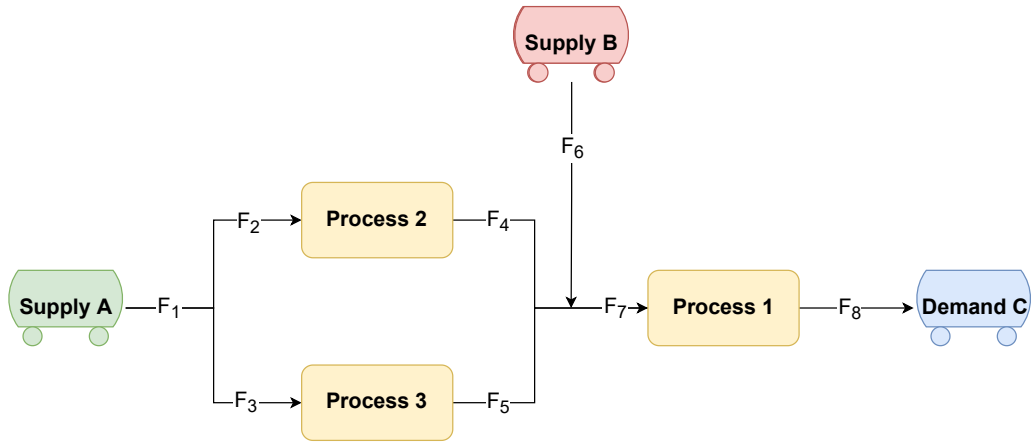

**Figure 4.3.1.** Process superstructure for Example 4.3.

The hierarchical GDP model is given as follows. The material balance constraints in the two stream junction points are given in (48) and (49), where $F_{s,t}$ is the flow (tons) in stream $s$ in period $t$ (where $t$ is in years). The amount of imported B and exported C are constrained by (50) and (51), respectively.

$$F_{1,t} = F_{2,t} + F_{3,t} \qquad \forall t \in T \quad (48)$$
$$F_{7,t} = F_{4,t} + F_{5,t} + F_{6,t} \qquad \forall t \in T \quad (49)$$
$$F_{6,t} \leq 5 \qquad \forall t \in T \quad (50)$$
$$F_{8,t} \leq 1 \qquad \forall t \in T \quad (51)$$

The installation and planning decisions are made in the nested disjunction given in (52), where the top-level decision is to install Process $i$ or not ($Y_i$ or $\neg Y_i$). If a process is installed, the respective nonlinear production yield constraint is enforced, where $g_1(F_{7,t}) = 0.9 \cdot F_{7,t}$, $g_2(F_{2,t}) = \ln(1 + F_{2,t})$, and $g_3(F_{3,t}) = 1.2 \cdot \ln(1 + F_{3,t})$. A process capacity balance is also applied to update the current capacity, $Q_{i,t}$, with the capacity in the previous period and the current capacity expansion, $QE_{i,t}$. The secondary level decision is to operate the installed process, $N_{i,t}$, or not. If the process is operated in period $t$, the exit flow is bounded by the process capacity, and the operating cost, $CO_{i,t}$, is determined with the parameter $\gamma_i$ ($\gamma_1 = \$900$, $\gamma_2 = \$1,000$, and $\gamma_3 = \$1,200$). The tertiary level decision is to expand the process capacity, $Z_{i,t}$, or not. The expansion cost, $CE_{i,t}$, is calculated with the fixed cost parameter, $\alpha_i$ ($\alpha_1 = \$3,500$, $\alpha_2 = \$1,000$, and $\alpha_3 = \$1,500$), and the variable cost parameter, $\beta_i$ ($\beta_1 = \$1,200/ton$, $\beta_2 = \$700/ton$, and $\beta_3 = \$1,100/ton$). It should be noted that each of the parameters used can also be indexed by time period if desired. Also, since the operating ($N_{i,t}$) and expansion ($Z_{i,t}$) decisions are nested, the second disjunct in each disjunction is not effectively a negation of the Boolean variable in the first disjunct (e.g., if process $i$ is not used, $\neg Y_i$, then both $N_{i,t} = False$ and $\neg N_{i,t} = False$ since these belong to the $Y_i$ disjunct). This is made clear in the cardinality clauses in (54) and (55). On the other hand, $\neg Y_i$ is effectively the negation of $Y_i$, as supported by (53).

$$\left[\begin{array}{c} Y_i \\ F_{s,t} = g_i(F_{s',t}) \; \forall t \in T \\ Q_{i,t} = Q_{i,t-1}|_{t>1} + QE_{i,t} \; \forall t \in T \\ \left[\begin{array}{c} N_{i,t} \\ F_{s,t} \leq Q_{i,t} \\ CO_{i,t} = \gamma_i \end{array}\right] \vee \left[\begin{array}{c} \neg N_{i,t} \\ F_{s,t} = 0 \\ CO_{i,t} = 0 \\ QE_{i,t} = 0 \\ CE_{i,t} = 0 \end{array}\right] \; \forall t \in T \\ \left[\begin{array}{c} Z_{i,t} \\ CE_{i,t} = \alpha_i + \beta_i \cdot QE_{i,t} \end{array}\right] \vee \left[\begin{array}{c} \neg Z_{i,t} \\ QE_{i,t} = 0 \\ CE_{i,t} = 0 \end{array}\right] \end{array}\right] \vee \left[\begin{array}{c} \neg Y_i \\ F_{s,t} = 0 \; \forall t \in T \\ F_{s',t} = 0 \; \forall t \in T \\ Q_{i,t}, QE_{i,t} = 0 \; \forall t \in T \\ CO_{i,t}, CE_{i,t} = 0 \; \forall t \in T \end{array}\right] \quad (52)$$

$$\forall i \in I, s \in S_i^{out}, s' \in S_i^{in}$$

$$\Xi(1, \{Y_i, \neg Y_i\}) \qquad \forall i \in I \quad (53)$$
$$\Xi(Y_i, \{N_{i,t}, \neg N_{i,t}\}) \qquad \forall i \in I, t \in T \quad (54)$$
$$\Xi(N_{i,t}, \{Z_{i,t}, \neg Z_{i,t}\}) \qquad \forall i \in I, t \in T \quad (55)$$

Additional logic constraints are given in (56)-(59). The cardinality clause in (56) allows installing *at most 1* of Process 2 or Process 3. This is equivalent to the proposition $\neg Y_2 \vee \neg Y_3$ used in the original paper, but generalizes for cases in which there are more than two potential processes in parallel. The implication in (57) ensures that Process 1 is installed if either Process 2 or Process 3 are installed. (58) and (59) enforce that process $i$ operate at least once if installed, with at least one expansion event scheduled between the beginning of the planning horizon (period 1) and each period $t$ in which the process is operated, respectively.

$$\Gamma(1, \{Y_2, Y_3\}) \qquad (56)$$
$$Y_i \Rightarrow Y_1 \qquad \forall i \in \{2,3\} \quad (57)$$
$$Y_i \Rightarrow \bigvee_{t \in T} N_{i,t} \qquad \forall i \in I \quad (58)$$
$$N_{i,t} \Rightarrow \bigvee_{t'=1}^{t} W_{i,t'} \qquad \forall i \in I \quad (59)$$

The variable domains are given in (60)-(66), where $F_s^{UB} = 5 \; ton \; \forall s \in S$, $QE_1^{UB} = 0.4 \; ton$, $QE_2^{UB} = 0.3 \; ton$, and $QE_3^{UB} = 0.3 \; ton$.

$$0 \leq F_{s,t} \leq F_s^{UB} \qquad \forall s \in S, t \in T \quad (60)$$
$$0 \leq Q_{i,t} \leq QE_i^{UB} \cdot t \qquad \forall i \in I, t \in T \quad (61)$$
$$0 \leq QE_{i,t} \leq QE_i^{UB} \qquad \forall i \in I, t \in T \quad (62)$$
$$0 \leq CE_{i,t} \leq \alpha_i + \beta_i \cdot QE_i^{UB} \qquad \forall i \in I, t \in T \quad (63)$$
$$0 \leq CO_{i,t} \leq \gamma_i \qquad \forall i \in I, t \in T \quad (64)$$
$$N_{i,t}, Z_{i,t} \in \{True, False\} \qquad \forall i \in I, t \in T \quad (65)$$
$$Y_i \in \{True, False\} \qquad \forall i \in I \quad (66)$$

The objective function is to minimize the system cost, as given in (67), where the stream costs, $p_s$, are given in **Table 4.3.1**. The model for Example 4.3 is thus given by (48)-(67).

$$\min \sum_{t \in T} \left( \sum_{s \in S} p_s \cdot F_{s,t} + \sum_{i \in I} (CO_{i,t} + CE_{i,t}) \right) \tag{67}$$

**Table 4.3.1.** Stream costs, $p_s$, in $/ton.

| Stream ($s$) | 1 | 2 | 3 | 4 | 5 | 6 | 7 | 8 |
|---|---|---|---|---|---|---|---|---|
| Cost ($/ton) | 1,800 | 0 | 0 | 300 | 100 | 7,000 | 0 | -10,800 |

There are some differences between this formulation and the one in the original paper by van den Heever and Grossmann (1999). The original formulation has the process capacity evolution constraint in the disjunct governed by $Z_{i,t}$. This requires specifying a new constraint, $Q_{i,t} = Q_{i,t-1}$, for the disjunct governed by $\neg Z_{i,t}$, which would also be required for the disjunct governed by $\neg N_{i,t}$. This is avoided by moving the process capacity balance to the upper-level constraints in $Y_i$. The same is true for the yield constraint, which we move from the $N_{i,t}$ disjunct to the $Y_i$ disjunct constraints. This requires that we only constrain the flow exiting the process in the secondary level disjunction, rather than both the entrance and exit flows. It is also more intuitive to specify the yield constraints when the processes are selected. Another major difference is that the original model does not use the cardinality constraints in (54) and (55). Instead, it uses the logic propositions (68) and (69). These propositions are contained in (54) and (55), but do not establish a proper hierarchical relationship since there is no link between $\neg N_{i,t}$ and $Y_i$, and $\neg Z_{i,t}$ and $N_{i,t}$.

$$N_{i,t} \Rightarrow Y_i \qquad \forall i \in I, t \in T \tag{68}$$
$$Z_{i,t} \Rightarrow N_{i,t} \qquad \forall i \in I, t \in T \tag{69}$$

An important thing to note is that the model in Example 4.3 is an example of a type of hierarchical GDP, that need not be hierarchical at all. This occurs when every disjunction has only two disjuncts, representing an *on* and an *off* state, where the *off* state has all relevant variables set to zero. When this occurs, (52) can actually be split into three sets of disjunctions without adding the "slack" disjunct observed in the Equivalent Single-Level GDP model. These three sets of disjunctions are given in (70)-(72). The cardinality constraints in (54)-(55) can be replaced by (68)-(69). The model composed of (48)-(51), (53), and (56)-(72) is referred to here as the **Non-hierarchical** formulation.

$$\begin{bmatrix} Y_i \\ F_{s,t} = g_i(F_{s',t}) \; \forall t \in T \\ Q_{i,t} = Q_{i,t-1}|_{t>1} + QE_{i,t} \; \forall t \in T \end{bmatrix} \vee \begin{bmatrix} \neg Y_i \\ F_{s,t} = 0 \; \forall t \in T \\ F_{s',t} = 0 \; \forall t \in T \\ Q_{i,t} = 0 \; \forall t \in T \\ QE_{i,t} = 0 \; \forall t \in T \end{bmatrix} \qquad \forall i \in I, s \in S_i^{out}, s' \in S_i^{in} \tag{70}$$

$$\begin{bmatrix} N_{i,t} \\ F_{s,t} \leq Q_{i,t} \\ CO_{i,t} = \gamma_i \end{bmatrix} \vee \begin{bmatrix} \neg N_{i,t} \\ F_{s,t} = 0 \\ CO_{i,t} = 0 \end{bmatrix} \qquad \forall i \in I, s \in S_i^{out}, t \in T \tag{71}$$

$$\left[ \begin{array}{c} Z_{i,t} \\ CE_{i,t} = \alpha_i + \beta_i \cdot QE_{i,t} \end{array} \right] \vee \left[ \begin{array}{c} \neg Z_{i,t} \\ QE_{i,t} = 0 \\ CE_{i,t} = 0 \end{array} \right] \qquad \forall i \in I, t \in T \qquad (72)$$

The Nested GDP model is compared against its Equivalent Single-Level formulation, and the Non-hierarchical formulation, by reformulating each of these into mixed-integer nonlinear programs (MINLPs) using the Hull reformulation, and solving them with in three modes: 1) relaxing the integrality constraints (NLP relaxation), 2) solving the MINLP with local search and range reduction disabled in BARON, and 3) solving the MINLP with local search and range reduction enabled. Since the models are nonlinear, the perspective functions were reformulated using the $\epsilon$-approximation from Furman et al. (2020), with $\epsilon = 10^{-6}$, which is the default reformulation approach in the *DisjunctiveProgramming* code. The model statistics are given in **Table 4.3.2**, where it is seen that the Nested formulation has a CPU time that is one order of magnitude smaller than that of the Equivalent Single-Level formulation when range reduction and local search are disabled. When range reduction and local search are switched on, the Nested formulation is faster than the Equivalent Single-Level formulation by only a factor of 2.5. The continuous relaxations for the Nested and Non-hierarchical formulations are equal (29% gap) and tighter than that of the Equivalent Single-Level formulation (69% gap). The performance of the Nested formulation is comparable to that of the Non-hierarchical one, with the latter having less continuous variables and constraints, and taking less time to solve when local search and range reduction are enabled. It is clear that for this model structure, a hierarchical formulation is not required.

**Table 4.3.2.** Model sizes and computational results of the MINLP models resulting from the Hull reformulations of the Equivalent Single-Level, Nested, and Non-hierarchical GDP models.

|  | Hull Reformulation | | |
| --- | --- | --- | --- |
|  | Equivalent | Nested | Alternate |
| **Model Size** | | | |
| *Binary Variables* | 384 | 258 | 258 |
| *Continuous Variables* | 2,499 | 2,058 | 1,554 |
| *Constraints* | 12,006 | 10,431 | 7,596 |
| **NLP Relaxation** | | | |
| *Relaxation Solution* | -$161,458 | -$122,829 | -$122,829 |
| *Relaxation Gap* | 69% | 29% | 29% |
| **MIP Solution**[a] | | | |
| *MIP Solution* | -$95,373 | -$95,373 | -$95,373 |
| *MIP Gap* | 0.0% | 0.0% | 0.0% |
| *BaR Iterations* | 267 | 7 | 7 |
| *Cuts Applied* | 69,831 | 1,898 | 1,912 |
| *CPU Time (s)* | 56.9 | 3.8 | 3.9 |
| **MIP Solution** | | | |
| *MIP Solution* | -$95,373 | -$95,373 | -$95,373 |
| *MIP Gap* | 0.0% | 0.0% | 0.0% |
| *BaR Iterations* | 1 | 1 | 1 |
| *Cuts Applied* | 946 | 950 | 965 |
| *CPU Time (s)* | 8.5 | 3.5 | 2.3 |

The optimal expansion profile is given in **Figure 4.3.2**, where it can be seen that Process 2 is not installed, but Processes 1 and 3 are, where the capacity in Process 1 increases to 1 ton/year by the third year, and Process 3 increases to 1.11 ton/year by the fourth year. The optimal system cost is -$95 thousand, meaning that plant generates profit.

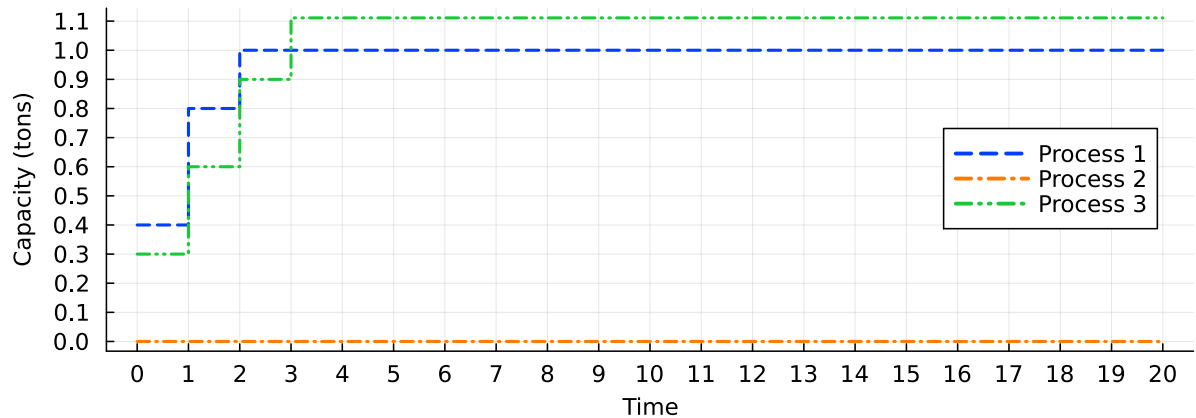

**Figure 4.3.2.** Capacity expansion profiles for each of the processes in Example 4.3.

## 5. Conclusions

Two main contributions are made in this paper to the generalized disjunctive programming (GDP) modeling framework. The first one is to add cardinality rules to the logic constraints to allow for constraints of the form *choose exactly m* Boolean variables to be *True* (or *at least m*, or *at most m*). For more than two Boolean variables, modeling these types of constraints via propositional logic (zeroth-order logic) is either cumbersome or not possible. Thus, introducing predicate logic (first-order logic) to express this new constraint form in GDP adds more expressiveness to logic-based models. The second contribution is to extend GDP for modeling hierarchical systems via nested disjunctions. Such an approach results in more intuitive models, but had not been formalized in the past, as classical GDP does not consider disjunction nesting. The notation and logic constraints for such structures is provided, along with theoretical proofs to the tightness of such models, versus equivalent single-level GDP models. It is shown that mixed-integer programming reformulations of nested GDP models have tighter relaxations than the reformulations of their single-level counterparts in both the Hull reformulation, as well as the Big-M reformulation when tight *M* values are used. However, when large *M* values are used, the reformulated nested models show worse performance due to the presence of multiple large *M* parameters in the nested constraints. Finding tight *M* values requires additional work, and can be done by applying interval arithmetic when the models are linear. However, for nonlinear models, a separate optimization model must be solved for each constraint to find the tightest *M* values. A discussion on using basic steps to flatten nested GDP models is also given, where flattening nested structures via basic steps improves performance in the Hull reformulations when few disjunctions are nested, but quickly becomes intractable as the number of nested disjunctions and disjuncts increases.

Three examples are presented to show the advantages of using nested structures. In the first example, the tightness of the continuous relaxations of nested linear models are compared geometrically with the relaxations of equivalent single-level models. The relaxations of models that preserve nested structures result in smaller feasible regions than their single-level counterparts, generally yielding significant computational savings. Example 4.2, a linear GDP, and Example 4.3, a nonlinear GDP, illustrate the computational advantages of nested GDP models for problems that integrate superstructure design, technology selection, and operations scheduling, and superstructure design, long-term operations planning, and capacity expansion planning, respectively. It is also shown that for systems with bi-disjunct constraints (disjunctions with only two disjunctions), where one disjunct represents an *off* state with all pertinent variables set to zero (e.g., zero flow), there is no advantage to modeling such systems as hierarchical, even when there may be several levels of decisions. Such systems can be modelled more simply with single-level disjunctions and the necessary linking constraints.

Future work includes investigating how explicit hierarchical structures can be exploited for informed model decomposition methods and branching strategies. Exploring applications of hierarchical GDP to other fields, such as decision trees and stochastic optimization with event constraints, is another potential area for development.

### Acknowledgment

The authors gratefully acknowledge the financial support from the Center of Advanced Process Decision-making at Carnegie Mellon University.

### Supplementary Material

All source code for the figures and examples in this paper can be found at https://github.com/hdavid16/Extensions-to-GDP-paper (repository will be made publicly available after manuscript acceptance).

### Nomenclature

The symbols for the sets, parameters, and variables used in Examples 4.2 and 4.3 are described below,

|  | Description |
|---|---|
| **Sets** |  |
| $i \in I$ | Processes |
| $j \in J$ | Technologies |
| $k \in K$ | Resources |
| $k \in K^{react}$ | Reactors |
| $k \in K^{tank}$ | Tanks |
| $s \in S$ | Streams |
| $s \in S_x^{in}$ | Streams entering $x$ |
| $s \in S_x^{out}$ | Streams exiting $x$ |
| $t \in T$ | Time periods |
| $u \in U_k$ | Number of installed units for resource $k$ |
| **Parameters** |  |
| $\alpha$ | Fixed installation/expansion cost |
| $\beta$ | Variable installation/expansion cost coefficient |
| $\gamma$ | Variable operating cost coefficient |
| $\nu$ | Yield coefficient |
| $\tau_i$ | Processing time for process $i$ |
| $F_s^{UB}$ | Upper bound on stream $s$ |
| $p_x$ | Price/cost of stream/material $x$ |
| $Q_k^{UB}$ | Upper bound on the capacity of resource $k$ |
| $QE_i^{UB}$ | Upper bound on the capacity expansion of process $i$ |

**Continuous Variables**

| | |
|---|---|
| $B_{i,t}$ | Total batch size for process $i$ starting in period $t$ |
| $\hat{B}_{i,t}$ | Unit batch size for process $i$ starting in period $t$ |
| $CE_{i,t}$ | Expansion cost for process $i$ in period $t$ |
| $CI_k$ | Installation cost for resource $k$ |
| $CO_{i,t}$ | Operating cost of process $i$ in period $t$ |
| $F_{s,t}$ | Flow in stream $s$ in period $t$ |
| $L_{k,t}$ | Level in tank $k$ in period $t$ |
| $\hat{L}_k$ | Level slack in tank $k$ in the final period (end of scheduling horizon) |
| $Q_k$ | Capacity of resource $k$ |
| $Q_{i,t}$ | Capacity of process $i$ in period $t$ |
| $QE_{i,t}$ | Capacity expansion for process $i$ in period $t$ |
| $R_{k,t}$ | Availability of resource $k$ in period $t$ |
| $\Delta R_{i,k,t}$ | Number of resources of type $k$ consumed for process $i$ in period $t$ |

**Boolean Variables**

| | |
|---|---|
| $N_{i,k,t,u}$ | Process $i$ is started on $u$ units of resource $k$ at period $t$ |
| $N_{i,t}$ | Process $i$ is operated in period $t$ |
| $\widehat{W}_{i,j}$ | Technology $j$ is used for process $i$ |
| $W_{i,k}$ | Resource $k$ is assigned to process $i$ |
| $X_{k,u}$ | Resource $k$ installed with $u$ units |
| $Y_i$ | Process $i$ is installed |
| $Z_{i,t}$ | Process $i$ undergoes a capacity expansion in period $t$ |